\documentclass[11pt]{article}
\usepackage[dvips]{graphicx}
\begin{document}
\tolerance=10000

%%\null

\hyphenpenalty=2000
\hyphenation{visco-elastic visco-elasticity}
\setcounter{page}{1}
\thispagestyle{empty}

%%\hoffset= 0.25truecm  %% TO RIGHT ->

%%%% FONTS

% ------------- specific fonts

\font\note=cmr10 at 10 truept  %% FOR ABSTRACT

%%%%%%%%%%%%%%%%%%%%%%%%%%%%%%%%%%%%%%%%%%%%%%%%
\newcommand{\eproof}{\rule{0.2cm}{0.2cm}}
%%\noindent
%%%%%%%% DEFINITIONS by Yu LUCHKO
\newcommand{\FTS}[2]{\frac{{\textstyle #1}}{{\textstyle #2}}}
 \newcommand{\NN}{\bf N}
 \newcommand{\ZZ}{\bf Z}
 \newcommand{\CC}{\bf C}
 \newcommand{\RR}{\bf R}
 \newcommand{\intl}{\int\limits}
\newcommand{\suml}{\sum\limits}
\setcounter{page}{1}
\thispagestyle{empty}
%%%%  DEFINITIONS by F. MAINARDI  %%%%%%%%%%%%%%%%%%%%%%%%%%%%%%
%% \font\title=cmbx10 scaled \magstep2
\font\bfs=cmbx10 scaled \magstep2
\def\eg{{\it e.g.}\ } \def\ie{{\it i.e.}\ }
\def\sg{\hbox{sign}\,}
\def\sgn{\hbox{sign}\,}
\def\sign{\hbox{sign}\,}
\def\e{\hbox{e}}
\def\exp{\hbox{exp}}
\def\ds{\displaystyle}
\def\dis{\displaystyle}
\def\q{\quad}    \def\qq{\qquad}
\def\lan{\langle}\def\ran{\rangle}
\def\l{\left} \def\r{\right}
\def\lra{\Longleftrightarrow}
\def\arg{\hbox{\rm arg}}
\def\d{\partial}
 \def\dr{\partial r}  \def\dt{\partial t}
\def\dx{\partial x}   \def\dy{\partial y}  \def\dz{\partial z}
\def\rec#1{{1\over{#1}}}
\def\log{\hbox{\rm log}\,}
\def\erf{\hbox{\rm erf}\,}     \def\erfc{\hbox{\rm erfc}\,}
%%%%%%%%%%%%%%%%%%%%%%%%%%%%%%%%%%%%55
%%%%%%%%%%%%%%%%%%%%%%%%%%%%%%%%%%%%%%%%%%%%%%%%%%%%%%
\def\L{{\cal L}} %%% Laplace Transform !!!!
\def\F{{\cal F}} %%% Fourier Transform !!!!
\def\M{{\cal M}}  %%% Mellin Transform
\def\I{{\cal I}}  %% The generic open interval  \in \RR
%%%%%%%%%%%%%%%%%%%%%%%%%%
\def\pni{\par \noindent}
\def\vsh{\smallskip}
\def\vs{\medskip}
\def\vvs{\bigskip}
\def\vvvs{\bigskip\medskip} %% {\vskip 1.5truecm}
\def\vsp{\par\noindent}
\def\vsn{\vsh\pni}
\def\cen{\centerline}
\def\ra{\item{a)\ }} \def\rb{\item{b)\ }}   \def\rc{\item{c)\ }}
%%%%%%%%%%%%%%%%%%%
%% DEF from PAMPLONA
\newcommand{\lt}{\left}
\newcommand{\rt}{\right}
\newcommand{\R}{{\mathcal R}}
\newcommand{\Mt}{\,\stackrel{\mathcal M}{\longleftrightarrow}\,}
\newcommand{\Ft}{\,\stackrel{\mathcal F}{\longleftrightarrow}\,}
\newcommand{\Lt}{\,\stackrel{\mathcal L}{\longleftrightarrow}\,}
%%%%%%%%
\def\ZZ{{\bf Z}}
%%% PROBABILITY,  EXPECTATION Value and VARIANCE
\def\PP{{\rm I\hskip-2pt P}\,}
\def\EE{{\rm I\hskip-2pt E}\,}
\def\Var{\hbox{Var}}
%%%%%%% making running heads on even and odd pages %%%%%%%
 \cen{FRACALMO PRE-PRINT   {\bf www.fracalmo.org}}
\vsh
\cen{{\bf Fractional Calculus and Applied Analysis}
  Vol. 8 No 1 (2005), pp. 7-38}
\vsh
\cen{www.diogenes.bg/fcaa/ \ ISSN \ 1311-0454}
\vs
\hrule
%%%%%%%%%%%%%%%%%%%%%%%%%%%%%%%%%%%%%%%%%%%%%%%%%%%%%%%%%%%%%%%%%%%%%%%%%
%%BEGINNING OF TEXT
%%%%%%%%%%%%%%%%%%%%%%%%%%%%%%%%%%%%%%%%%%%%%%%%%%%%%%%%%%%%%%%%%%%%%%%%%
   \vskip 0.50truecm

%%%%%%%%%%%%  TITLE AND ABSTRACT %%%%%%

\font\title=cmbx12 scaled\magstep2
\font\bfs=cmbx12 scaled\magstep1

%% \vsh
\begin{center}
%%  {\bf First International Conference on
%%  \\  "Nonlinear Analysis and Mechanics of Continuous Media"
%% \\  Ho Chi Minh City, Vietnam, 7-10 January 2003.}\\

{\title Renewal processes}
\vs

{\title of Mittag-Leffler and Wright type}

\vvs

 {Francesco MAINARDI}$^{(1)}$,
{Rudolf GORENFLO}$^{(2)}$,
{Alessandro VIVOLI} $^{(1)}$

%%%%%%%%%%%%%%%%

\vs

$\null^{(1)}$
 Dipartimento di Fisica, Universit\`a di Bologna and INFN, \\
              %%%  Sezione di Bologna, \\
Via Irnerio 46, I-40126 Bologna, Italy \\
%% Tel: +39-051.2091098 $\;$ Fax: +39-051.247244 $\;$\\
E-mail: {\tt mainardi@bo.infn.it} $\;$ URL: {\tt www.fracalmo.org}
\\ [0.25 truecm]

$\null^{(2)}$
 First Mathematical Institute,
 %%%Department of Mathematics and  Informatics,
 Free   University of Berlin, \\
 Arnimallee  3, D-14195 Berlin, Germany \\
%% Tel: +49-30-838.75.352 $\;$ Fax: +49-30-838.75.403 $\;$ \\
E-mail: {\tt gorenflo@mi.fu-berlin.de}
\\ [0.25 truecm]

%%%%%%% DATE to be omitted in the FINAL VERSION
%% \\ [0.25truecm]
 \def\date#1{\gdef\@date{#1}} \def\@date{\today}
%% {\bf Ver Final, \@date }
%% {\bf REVISED VERSION:  \@date}
%%%%%%%%%%%%%%%%%%%%%%

\end{center}

\cen{\bf Abstract} %%% Abstract follows

\vskip 0.1truecm
\noindent

 After sketching the basic principles of renewal theory we
discuss the classical Poisson process and offer two other processes,
namely the renewal process of Mittag-Leffler type and the renewal process
of Wright type, so named by us because special functions of
Mittag-Leffler and of Wright type appear in the definition of the relevant
waiting times. We compare these three processes with each other,
furthermore consider corresponding renewal processes with reward and
numerically their long-time behaviour.

\vskip 0.20truecm
\noindent
{\it 2000 MSC}: %%%%   Mathematics Subject Classification}:
  26A33, 33E12, 33E20, 44A10, %% 44A20,
  44A35, %% 45E10, 45J05, 45K05,  47G30,
  60G50, 60J05, 60K05.
  %% 60E07, 60J60, 60J65, 73F05.
%%%%%%%%%
%% 33E12 is the new AMS 2000 Classification for the M-L function
%% 60G52 is the AMS class for stable process, suggested by GORENFLO!
%\begin{ams}  %% 60E20, 49G03, 49F10
%26A33,  %%%%  (main);    Fractional derivatives and integrals
%% 33E12, %% Mittag-Leffler type functions
%%  33C60,  %% hypergeometric integrals and functions defined by them
%% 44A10,  %% Laplace Transforms
% 45K05,  %% integro-partial differential equations
%% 26A33, %% Fractional derivatives and integrals
%%%% 33B20, %% Incompleta beta and gamma functions (error functions,..)
%% 33E12, %% Mittag-Leffler type functions
%%%% 33C60,  %% hypergeometric integrals and functions defined by them
%%%% 42A38   %% Fourier and Fourier-Stieltjes transforms
%% 44A10, %% Laplace Transforms
%% 44A15, %% Special transforms (Legendre, Hilbert, etc.)
%% 44A35, %% Convolution
%%%% 45K05,  %% integro-partial differential equations
%%%% 47G30,  %% Pseudo-differential operators
%%%% 60-xx  Probability theory and  Stochastic processes !!!!!!!!
%%%% 60G18, %%%  Self-similar processes
%% 60G50, %% sums of independent random variables, random walks
%%%% 60G51, %% processes with independent increments (Levy processes!)
%%%% 60G52, %% stable processes
%% 60G55, %% Stable processes (?), Point processes
%% 60J05, %% Markov processes with discrete parameters
%%%% 60J25  %%    Markov processes with continuous parameters
%%%% 60J60. %% diffusion processes
%%%%  60J70 applications of diffusion processes
%% 60K05. %% Renewal theory
%% 60K25 %% Queueing theory
\vskip 0.25truecm
\noindent
{\it This paper is dedicated to Acad. Bogoljub Stankovi{\'c},  Emeritus Professor,
Department  of Mathematics, University of Novi Sad, Serbia,%%  21000 Novi Sad, Serbia,
on the occasion of his 80-th birthday (November 1,  2004)}
%%%%%%%%%% \vfill\eject

%%\section{Introduction}
\section*{1. Introduction}

%% The purpose of this  paper  is to point out
%% the relevance  of functions of  Mittag-Leffler and Wright type
%% in renewal theory for constructing   subordinated stochastic processes
%% of fractional diffusion.
%% A) We provide a  generalization of and a variant to the Poisson
%%   process (without and with reward) which is known to play a
%%   fundamental role in renewal theory. "
%% B) It is well known that the Poisson process (with and without reward)
%% plays a fundamental role in renewal theory. Here we provide
%% a generalization and a variant to this process. "
It is well known that the Poisson process (with and without reward)
plays a fundamental role in renewal theory.
In this paper, by means of functions of  Mittag-Leffler and Wright type
  we provide a  generalization of and a variant to this classical process
and construct interesting subordinated stochastic processes
of fractional diffusion.
%%%%%%%%%%%%%%%%%\vfill\eject

The  plan of the paper is as  follows.
%%%%%%%%%% \vfill \eject
%%%%%%%%%%%%%%
In Section 2 we recall the basic renewal theory including its
fundamental concepts like waiting time between events, the survival
probability, the renewal function.
%%%%%%%%%%

If the waiting time is exponentially distributed we  have the classical
Poisson process, which is Markovian: this is the topic of Section 3.
However, other waiting time distributions are also relevant in
applications, in particular such ones with a fat tail caused by a power
law decay of their density.
In this context we analyze, respectively in Sections 4 and 5,
  two non-Markovian
renewal processes with  waiting time distributions described by
functions of Mittag-Leffler and Wright type, that exhibit a similar
power law decay.
They both depend on a parameter $\beta \in (0, 1)$ related to the
common exponent in the power law.
In the limit $\beta = 1$ the first becomes the Poisson process
whereas the second goes over into the deterministic process
producing its events at equidistant instants of time.

In Section 6, after sketching the basic differences between
the renewal processes of Mittag-Leffler and Wright type,
 we compare numerically their survival functions and their
probability densities in the special case $\beta =1/2$
with respect to  the corresponding functions of
the classical Poisson process.

In Section 7
 we discuss the general renewal process with reward
(also called the {\it compound renewal process}), interpret it
physically as a {\it continuous time random walk} (CTRW),
and review  how the resulting sojourn probability density function,
evolving  in time according to an integral equation,
can analytically be represented as an infinite series.
We then  consider the {\it compound} renewal processes
of Mittag-Leffler and Wright type including the limiting case $\beta =1$.
%% These processes, because of their asymptotic properties
%%  are known  to reduce, in the diffusion limit,
%% to  {\it space-time  fractional diffusion processes}.
%% whose spatial probability distributions are governed
%% by generalized diffusion equations with partial
%% derivative of fractional order in space and time.
%%%%%%%%
In Section 8, we recall the conditions
of the well-scaled transition to the diffusion limit
according to which the  CTRW integral equation
of our compound processes
 reduces to  the time fractional diffusion equation (TFDE).
For  all three cases (Poisson, Mittag-Leffler, Wright),
by numerically summing the corresponding series for
the cumulative probability function, we produce plots of its behaviour in
time. These plots nicely illustrate that with increasing time
the behaviour of time-fractional diffusion is approached.
%%%%%%%%%
%% When, in particular, the stochastic process of the reward
%% is spatially-symmetric with finite variance,
%% the integral equation  of the corresponding CTRW
%% Taking  the same probability
%% distributions for the waiting time considered in Section 6,
%% we compare the spatial probability distributions evolving in time
%% provided by the  fundamental solution
%% of the integral equation of the CTRW
%% with the fundamental solutions of the corresponding limiting  TFDE.
%%
Finally, concluding remarks are given in Section 9.
%%%%%%%%\vfill\eject
%%%%%%%%%%

%% \section{\bfs Essentials of  renewal theory}
\section*{2. Essentials of  renewal theory}

For the sake of reader's convenience,
we present a brief  introduction  to the renewal theory
by using our notation.
For more details see \eg the classical
treatises by
%%  Khintchine \cite{Khintchine QUEUEING60},
Cox \cite{Cox RENEWAL67},
%% Gnedenko \& Kovalenko \cite{Gnedenko-Kovalenko QUEUEING68},
Feller \cite{Feller 71},
and the more  recent book by Ross \cite{Ross PROBMOD97}.
%%  and by Beichelt \& Fatti \cite{Beichelt-Fatti STOCHPROC02}.
%%%%%%%\vfill\eject
%%%%%%%%%%%
We begin to recall that a stochastic process $\{N(t), \; t\ge0\}$ is said
to be a {\it counting process}
if $N(t)$ represents the total number of "events"
that have occurred up to time $t$.
The concept of  {\it renewal process} has been developed as a
stochastic model for describing the class of counting processes
for which the times between successive  events,
$T_1,  T_2, \dots$, are
independent  identically distributed  ($iid$)
non-negative random variables, obeying a given probability law.
We call these times
 {\it waiting times} (or inter-arrival times)
and  the times
$$t_0 = 0\,, \q   t_k= \sum_{j=1}^k T_j\,, \q k \ge 1\,. \eqno(2.1)
    $$
 {\it renewal times}.
That is,
$t_0=0$ is the starting time,
$t_1 =T_1$ is the time of the first renewal,
$t_2 = T_1 +T_2$ is the time of the second renewal, and so on;
in general $t_k$ denotes the time of the $k$th renewal.
%%%%%%%%%%%

Let the {\it waiting times} be distributed like $T$ and let
$$ \Phi(t) := P (T \le t ) \eqno(2.2)$$
be the common  probability distribution function.
Here $P$ stands for {\it probability}.
We assume $\Phi(t)$ to be  absolutely continuous,
so that we can define its  probability density function
 $\phi(t)$ as
$$ \phi (t) = \frac{d}{dt} \Phi(t) \,,\q
   \Phi(t)  = \int_0^t \phi (t')\, dt'\,.
\eqno(2.3)$$
We recall that $\phi (t) \ge 0$ with $\int_0^\infty \phi(t)\,dt =1$
and $\Phi(t)$ is a non-decreasing function in $\RR^+$
with $\Phi(0) =0$,  $\Phi(+\infty) =1$.

Let us remark that, %%  for ease of language,
as it is popular in Physics, we use
the word density also for generalized functions
%%% in the sense of Gel'fand \& Shilov \cite{Gelfand-Shilov 64},
that can be interpreted as probability measures.
In these cases the function $\Phi(t)$ may lose its   absolute
continuity.
Often, especially in Physics,  the {\it probability density function}
is  abbreviated    by  $pdf$, so that, in order to avoid confusion,
the probability distribution function is called
the {\it probability cumulative  function} (being
the integral of the density) and abbreviated by $pcf$.
%%%%%%\vfill\eject
%%%%%%%%%%%%%
 When the nonnegative random variable  represents
 the lifetime of technical systems, it is common to
 call $\Phi(t)$  the {\it failure probability}
and
$$ \Psi(t) := P \l(T > t\r) = \int_t^\infty \phi (t')\, dt' = 1-\Phi(t)\,,
\eqno(2.4)$$
 the {\it survival probability}, because $\Phi(t)$ and $\Psi(t)$ are
the respective probabilities that the system does or does not fail
in $(0, t]$. These terms, however,  are commonly adopted for  any
renewal process.
%%%%%%%%%%

As a matter of fact the {\it renewal process} is defined by the
{\it counting process}
$$  N(t):=  \cases{
 0    & {for} $\; 0 \le t< t_1\,,$  \cr
   \hbox{max} \l\{k | t_k \le t, \;k =  1, 2, \dots\r\}
  & {for} $ \; t \ge t_1\,.$\cr
  }\eqno(2.5)$$
$N(t)$ is thus the random number of renewals
occurring in $(0,t]$. We easily recognize that
$\Psi(t) = P \l(N(t) =0\r)\,.$
%%%%%%%%%\vfill\eject
%%%%%%%%%
For an example of a renewal process, suppose
that  we have an infinite supply of light-bulbs
whose lifetimes are $iid$ %%% independent, identically distributed
random variables.
Suppose also that we use a single light-bulb at a time,
and when it fails we {\it immediately} replace it with a new one.
Under these conditions, $\{N(t), \, t\ge 0\}$ is a renewal process
when $N(t)$ represents the number of light-bulbs that have failed
by time $t$.
%%%%%%%%%%%%

%% For ease of language we also refer to $N(t)$ as
%% the {\it counting function}

Continuing in the general theory, we set
$F_1(t) = \Phi(t)$, $f_1(t) = \phi (t)$, and in general
$$ F_k(t) :=  P\l(t_k = T_1+ \dots +T_k \le t \r)\,,
 \; f_k(t) = \frac{d}{dt} F_k(t)\,, \; k\ge 1\,. \eqno(2.6)$$
$F_k(t)$ is the probability that the sum of the first $k$
waiting times does not exceed  $t$,
and $f_k(t)$ is the corresponding  density.
$F_k(t)$ is normalized
because
$         {\ds \lim_{t \to \infty}} F_k(t) =
  P\l(t_k = T_1+ \dots +T_k < \infty \r)= \Phi(+\infty) = 1.$
   In fact,    the sum of $k$ random variables each of which is
    finite with probability 1 is finite with probability 1 itself.
%%%%%

We set for consistency
 $F_0(t) = \Theta (t)$, the Heaviside unit step
function
(with $\Theta (0) := \Theta(0^+)$) so that $F_0(t) \equiv 1$ for $t\ge 0$,
 and $f_0(t) = \delta(t)$,
the Dirac delta  generalized function.
%%%%%%%%%%%%%%%%%%%%%%%%%%%%%%%%%%%%%%%%%%%%%%%%%%%%

A relevant quantity related to the counting process $N(t)$ is the
function $v_k(t)$ that represents the probability that
$k$ events occur in the closed interval $[0,t]$.
Using the basic assumption that the waiting times are $i.i.d.$
random variables, we get,
for any $k \ge 0$,
$$ v_k(t) :=P\l(N(t) =k\r) = P\l(t_k \le t \,,\,t_{k+1} > t\r)
=  \int_0^t f_k(t')\, \Psi (t-t')\, dt'\,.
\eqno(2.7)$$
We note that for $k=0$ we recover  $v_0(t) = \Psi(t)$.

%%%%%%%%%%%%%% RENEWAL FUNCTION %%%%%%%%%%%%%%%
Another  relevant quantity
 is the {\it renewal function} $m(t)$ defined as
the expected value of the process $N(t)$, that is
$$ m(t) := E \l(N(t)\r) = \langle N(t) \rangle =
   \sum_{k=1}^\infty P \l(t_k \le t\r)\,.
\eqno(2.8) $$
Thus this function
 represents the average number of events in the interval $(0, t]$
and  can be shown to uniquely determine the renewal process
\cite{Ross PROBMOD97}.
It is related
to the waiting time distribution by the so-called %% (famous)
{\it Renewal Equation},
$$ m(t) = \Phi(t) + \int_0^t m(t-t') \, \phi (t') \,dt' =
   \int_0^t \l[1+ m(t-t')\r] \, \phi (t') \,dt'\,.
\eqno(2.9)$$
%% Before treating a special distribution of waiting times, let us give a
%% few general properties. %% which are intuitive (but not so easy to prove).
If the mean waiting time (the first moment) is finite, namely
$$ \rho := \langle T \rangle  = \int_0^\infty t \, \phi (t) \, dt < \infty\,,
\eqno(2.10) $$
it is known that, with probability 1,
$ t_k/{k} \to \rho$  as $k \to \infty\,,$
and
$ N(t)/t \to 1/{\rho}$ as $t \to \infty\,.$
These facts imply the {\it elementary renewal theorem},
   $$\frac{m(t)}{t} \to \frac{1}{\rho} \,\q \hbox{as}\q t \to \infty\,.
\eqno(2.11)$$
%%%%%%%%%%%%%%%%%%%% \vfill\eject
%%%%%%%%%%
We shall also consider
 renewal processes
in which    the mean waiting time is infinite
because the waiting time laws have fat tails with
 power law asymptotics:
  $$\phi(t) \sim \frac{A_\infty}{t^{(1+\beta)}}\,, \;
        \Psi(t) \sim \frac{A_\infty}{\beta\, t^{\beta}}\,, \;
\hbox{for}
   \; t \to \infty\,, \q  0<\beta <1\,, \; A_\infty >0\,.
\eqno(2.12)$$

%%%%%%%%%%% NOTATION OF LAPLACE CONVOLUTION

It is convenient to use the common $\, *\,$ notation
for the Laplace convolution of two causal well-behaved
(generalized) functions $f(t)$ and $g(t)$,
$$   \int_0^t   f(t')\, g(t-t')\, dt' = \l(f \,*\, g\r) (t) =
     \l(g \,* \, f\r) (t)  = \int_0^t   f(t-t')\, g(t')\, dt'\,. $$
Being  $f_k(t)$ the $pdf$ of the sum of the  $k$
$iid$ random variables
$T_1,  \dots, T_k$
with $pdf$ $\phi (t)\,, $ we recognize that
$f_k(t)$ is the $k$-fold convolution of $\phi(t)$
with itself,
$$ f_k(t) =  \l(\phi^{*k}\r) (t)\,, \eqno(2.13)$$
so that Eq. (2.7)  simply reads:
$$  v_k(t) :=   P\l(N(t) =k\r) =
\l(\Psi \,*\, \phi^{*k}  \r)(t)\,. \eqno(2.14)$$
%%%%%%%%%%%
We  note that in the convolution notation the {\it renewal equation} (2.9)
reads
$$ m(t) = \Phi(t) + (m\,* \, \phi) (t)\,.\eqno(2.15) $$
%%%%%%%%%%%%%    LAPLACE TRANSFORMS %%%%%%%%%%%%%%%%%%%%%
The presence of Laplace convolutions allows us to treat
a renewal process by  the Laplace transform.

Throughout  we will denote
by   $\widetilde f(s)$
the Laplace transform
of a sufficiently well-behaved (generalized) function
$f(t)$  according to
$$   {\L} \l\{ f(t);s\r\}=   \widetilde f(s)
 = \int_0^{+\infty} \e^{\ds \, -st}\, f(t)\, dt\,,
\q s > s_0\,,
$$
and for $\delta (t)$
consistently we will have
$  \widetilde \delta (s) \equiv 1\,. $
{\it} Note that for our purposes we agree to take $s$ real.

We recognize  that
(2.13)-(2.14) reads in the Laplace domain
$$       \widetilde{f}_k(s) =  \l[{\widetilde \phi (s)} \r]^k \,, \q
  \widetilde{v}_k(s)
= \l[{\widetilde \phi (s)} \r]^k \,   \widetilde \Psi (s)
             \,,\eqno(2.16)$$
where, using (2.4),
$$ \widetilde \Psi (s)  =
\frac{ 1- \widetilde \phi (s)} {s}\,.
\eqno(2.17)$$
%%%%%%%%%%%%%%
Then, in the Laplace domain, the {\it renewal equation} reads
$$ \widetilde m(s) = \widetilde \Phi (s) +
  \widetilde m(s)\, \widetilde \phi (s) \,,\q
\hbox{with} \q
 \widetilde \Phi (s)= \frac{\widetilde \phi (s)}{s}\,, \eqno(2.18)$$
from which
$$ \widetilde m(s)  =
 \frac{\widetilde \phi (s)}{s\,\l[1-\widetilde \phi (s)\r]}
\,,
%% and  $$
\q  \widetilde \phi (s) =   \frac{s \, \widetilde m(s)}
         { 1 +  s\,\widetilde m(s)}  \,. \eqno(2.19)$$
%%%%%%%%%

%%%%%%%%%\vfill\eject%%%%%%%%%%

%% \section{{\bfs The Poisson process as a renewal process}}
\section*{3. The Poisson process as a renewal process}

The most celebrated   %%% simple
renewal process is the {\it Poisson process}
(with  parameter $\lambda >0$).
It is characterized by
 a survival function  of exponential type:
 $$ \Psi(t) = \e^{-\lambda t}\,, \q t \ge 0 \,.
\eqno(3.1)$$
As a consequence, the corresponding density for the waiting times
%% have the $pdf$
is exponential as well:
$$  \phi (t) = \lambda \, \e^{-\lambda t}\,, \q t\ge 0\,,
\eqno(3.2)$$
with   moments %%%
 $$ \langle T\rangle = \frac{1}{\lambda }\,, \q
   \langle T^2 \rangle = \frac{1}{\lambda^2 }\,,
\q \dots \,,  \q \langle T^n \rangle = \frac{1}{\lambda^n }\,, \q \dots \,,
\eqno(3.3)$$
and the  renewal function is  linear:
$$ m(t) = \lambda \, t\,,  \q t \ge 0
\,.\eqno(3.4)$$
The Poisson process is Markovian because
the exponential distribution is characteristic for
 a process without  memory.
%%%%%%%%%%%%%%%%%%%%%%%%%%%%%
We know that the probability  that $k$ events occur in an
interval of length $t $ is  the celebrated {\it Poisson distribution}
with parameter $\lambda t$,
$$ v_k(t):= P\l( N(t) = k\r) =   \frac{(\lambda t)^k}{k!} \, \e^{-\lambda t}
\,, \q t \ge 0\,, \q k = 0,1, 2, \dots\,. \eqno(3.5)$$
The probability distribution related to the sum of $k$ $iid$
exponential random variables  is known to be
the so-called {\it Erlang distribution}  (of order $k\ge 1$).
The corresponding density (the {\it Erlang} $pdf$) is thus
$$ f_k(t) = \lambda \,  \frac{(\lambda t)^{k-1}}{(k-1)!}
    \, \e^{-\lambda t}\,,\q t \ge 0 \,,
\q k =1,2, \dots \,, \eqno(3.6)$$
and the  corresponding {\it Erlang} $pcf$  is
$$ F_k(t) = \int _0^t f_k(t')\, dt' =
     1 -  \sum_{n=0}^{k-1}
      \frac{(\lambda t)^n}{n!} \, \e^{-\lambda t}  =
 \sum_{n=k}^{\infty}
      \frac{(\lambda t)^n}{n!} \, \e^{-\lambda t}\,,\q t \ge 0
\,. \eqno(3.7)$$
In the limiting case $k=0$ we recover
$f_0(t) = \delta(t)$
and  $F_0(t) = \Theta (t)$.
%%%%%%

The results  (3.4)-(3.7) can be easily obtained
by the Laplace transform technique. In fact, starting from
the Laplace transforms of the probability laws (3.1)-(3-2)
$$ \widetilde\Psi(s) =  \frac  {1 } {\lambda +s}\,, \q
 \widetilde\phi(s) =
\frac{\lambda } {\lambda +s}\,,
 \q     \eqno(3.8)$$
and, using (2.16)-(2.19),
we have,
$$  \widetilde m(s)  =
 \frac{\widetilde \phi (s)}{s\,\l[1-\widetilde \phi (s)\r]}
   = \frac{\lambda }{s^2}\,,\q
\widetilde v_k(s) =  [\widetilde\phi(s)]^k  \, \widetilde\Psi(s)
=  \frac{\lambda ^k}{(\lambda +s)^{k +1}}\,,
 \eqno(3.9)$$
hence (3.4)-(3.5), and
$$  \widetilde f_k(s) =  [\widetilde\phi(s)]^k =
    \frac  {\lambda^k } {(\lambda +s)^k}\,,
\q    \widetilde F_k(s) =  \frac{[\widetilde\phi(s)]^k}{s} =
    \frac  {\lambda^k } {s (\lambda +s)^k}\,,
\eqno(3.10)$$
hence (3.6)-(3.7).
%%%%%%%%%%%%%%%%%\vfill\eject
%%%%%%%%%%

%% \section{{\bfs The renewal  process of Mittag-Leffler type} }
\section*{4. The renewal  process of Mittag-Leffler type}

A "fractional" generalization of the   renewal Poisson process
has been recently proposed
by Mainardi, Gorenflo and Scalas \cite{Mainardi VIETNAM04}.
Noting that
 the survival probability
for the Poisson renewal process (with parameter $\lambda >0$)
obeys the ordinary differential equation  (of relaxation type)
$$     \frac{d}{dt} \Psi(t) = -\lambda \Psi(t)\,, \q t \ge 0\,; \q
\Psi(0^+) =1\,. \eqno(4.1)$$
the required generalization is obtained by
replacing in (4.1) the first derivative
by
%% the integro-differential operator $\, _tD_*^\beta$,
the fractional derivative
(in  Caputo's sense\footnote{%%
%%%%%%%%%%%%%% FOOTNOTE ON CAPUTO FRACTIONAL DERIVATIVE %%%%%
The Caputo derivative of order $\beta \in (0,1]$ of
 a well-behaved function $f(t)$ in $\RR^+$ is
$$
    _tD_*^\beta \,f(t) :=
\cases{
    {\ds \rec{\Gamma(1-\beta )}} {\ds\int_0^t}
 {\ds {f^{(1)}(\tau)\over (t-\tau )^{\beta }}\, d\tau}\,,
  & $\; 0<\beta  <1\,, $\cr\cr
     {\ds {d\over dt}} f(t)\,,
    & $\; \beta  =1\,. $\cr\cr }
    $$
Its Laplace transform turns out as
$$ {\cal L} \{ \, _tD_*^\beta  \,f(t) ;s\} =
      s^\beta \,  \widetilde f(s)
   -  s^{\beta  -1}\, f(0^+)\,.$$
For more information on  the theory and the applications
of the Caputo fractional derivative
(of any order $\beta >0$), see \eg
\cite{Caputo FCAA01,GorMai CISM97,Mainardi CISM97,Podlubny BOOK99}.}
%%%%%%%%%%%%%%%%%%%%%%%%%%%%%%%%%%%%%%%%%%%%%55
of order $\beta  \in (0,1]$.
We thus write, taking for simplicity $\lambda =1$,
$$      \, _tD_*^\beta \,\Psi(t) =
- \Psi(t)\,, \q t \ge 0\,,\q 0<\beta \le 1\,; \q
\Psi(0^+) =1\,. \eqno(4.2)$$

The solution  of Eq. (4.2) can be obtained by using the technique of
the Laplace transforms.
We have for $t\ge 0\,$:
$$\Psi(t) =  E_\beta (-t^\beta)\,, \q \hbox{from} \q
\widetilde\Psi(s) = \frac{s^{\beta -1}}{1+ s^\beta}\,,
 \q 0<\beta \le 1\,,
\eqno (4.3)$$
hence :
$$ \phi(t) =    -   \frac{d} {dt} \Psi(t) =
            -   \frac{d}{dt}  E_\beta (-t^\beta),
\; \hbox{corresponding to} \;
\widetilde\phi(s) = \frac{1}{1+ s^\beta}.
 \eqno (4.4)$$
%%%%%%%%%%%%%%%%%\vfill\eject
%%%%%%%%%%%%%%%
Here $ E_\beta$ denotes the Mittag-Leffler function\footnote{%%}
%%%%%%%%%%%%%%%%%   FOOTNOTE ON THE MITTAG-LEFFLER FUNCTION %%%
The Mittag-Leffler function with parameter  $\beta$
is defined as  %%%% in the complex plane by the power series
$$ E_\beta (z) :=
    \sum_{n=0}^{\infty}\,
   {z^{n}\over\Gamma(\beta\,n+1)}\,, \q \beta >0\,, \q z \in \CC\,.
 $$
%%%%
It is an entire function of order $\beta $
and reduces for $\beta=1$ to $\exp (z)\,.$
For detailed information on  the functions of Mittag-Leffler type
the reader  may  consult \eg
\cite{%% Djrbashian 66,
Erdelyi HTF,GorMai CISM97,Kiryakova 94,MaiGor JCAM00,%%%
Podlubny BOOK99,SKM 93} and references therein.} of order $\beta$.
%%%%%%%%%%%%%%%%%%%%%%%%%%%%%%
%%% In view of the presence of the Mittag-Leffler function in (4.3)-(4.4),
%% we agree to refer to    this  process as to the
We call this  process  the
{\it renewal process of Mittag-Leffler type} (of order $\beta $).

Hereafter, we find it convenient to summarize
the  most relevant features  of the functions $\Psi(t)$ and $\phi(t)$
when  $0< \beta <1\,.$
%% that turn out to be most relevant for our purposes.
We begin to quote their  expansions
in  power series of $t^\beta $ (convergent  for $t\ge 0$)
and  their asymptotic representations for  $t\to \infty $,
$$ \Psi(t)
   = {\ds \sum_{n=0}^{\infty}}\,
  (-1)^n {\ds {t^{\beta n}\over\Gamma(\beta\,n+1)}}
 \,\sim \,  {\ds {\sin \,(\beta \pi)\over \pi}}
  \,{\ds  {\Gamma(\beta)\over t^\beta}}\,,
 %%%% & $\, t\to \infty \,,$} \q 0<\beta <1\,,
     \eqno(4.5) $$
%% and
$$ \phi(t)
= {\ds {1\over t^{1-\beta}}}\, {\ds \sum_{n=0}^{\infty}}\,
  (-1)^n {\ds {t^{\beta n}\over\Gamma(\beta\,n+\beta )}}
 \, \sim \,  {\ds {\sin \,(\beta \pi)\over \pi}}
  \,{\ds  {\Gamma(\beta+1)\over t^{\beta+1}}}\,.
%%%  & $\, t\to \infty \,.$} \q 0<\beta <1\,,
     \eqno(4.6) $$
%%%%%%%%
%% The  expression for
%% $\phi(t)$ can be shown to be equivalent to that one
%% obtained in   \cite{Hilfer-Anton 95} in terms of the
%% generalized Mittag-Leffler function in two parameters.
%%%
We recognize that for $0<\beta <1$ both  functions
 $\Psi(t)$, $\phi(t)$
%% even if they lost their exponential decay by exhibiting power-law tails
%% for large times,
keep   the complete monotonicity\footnote{
%% of the exponential function of $\exp (-t)\,.$
%%%%%    FOOTNOTE ON COMPLETE MONOTONICITY OF MITTAG-LEFFLER
Complete monotonicity of  a function
 $f(t)$   means,  for $n=0,1,2,\dots$, and  $t\ge0$,
$ {\ds (-1)^n {d^n\over dt^n}\, f(t) \ge 0}$,
or equivalently, its representability as (real) Laplace transform
of a non-negative  %% (ordinary or generalized)
function or measure.
Recalling the theory of the Mittag-Leffler functions
of order less than 1,  we obtain
for $0<\beta <1$ the following representations, see \eg \cite{GorMai CISM97},
$$ \Psi(t)  =
  {\ds{\sin \,(\beta \pi)\over \pi}\,
   \int_0^\infty \!
   { r^{\beta  -1}\, \e^{\,\ds -rt}\over
    r^{2\beta } + 2\, r^{\beta } \, \cos(\beta  \pi) +1}\, dr}\,,
  \q t \ge 0\,,$$
%% and
$$  \phi(t) =
  {\ds {\sin \,(\beta \pi)\over \pi}\,
   \int_0^\infty \!
   { r^{\beta}\, \e^{\,\ds -rt} \over
    r^{2\beta } + 2\, r^{\beta}\,\cos(\beta \pi) +1}\, dr}\,,
 \q t\ge 0  \,.     $$},
 characteristic
for the Poissonian case $\beta=1$.
%%%%%%%%%%%%%%%%%%%%% \vfill\eject

%%%%%%%%%%%%%%%%%%
In contrast to the Poissonian case,
in the case  $0<\beta <1$ the mean waiting time
is infinite because the waiting time laws no longer decay exponentially
but exhibit power
law asymptotics according to (2.12)
where the constant $A_\infty$ derived from (4.5)-(4.6) is
$$     A_\infty =\Gamma(\beta+1)\, \sin (\beta \pi)/ \pi =
     \beta \,\Gamma(\beta)\, \sin (\beta \pi)/ \pi \,. \eqno(4.7)$$
%% %%%({\it power-law decay}).
As a consequence
 the process turns out to be no longer
Markovian but of long-memory type.
%%%%%%%%%%%%%%%%%%%
\begin{figure}%% [htbp]
 \includegraphics[width=.52\textwidth]{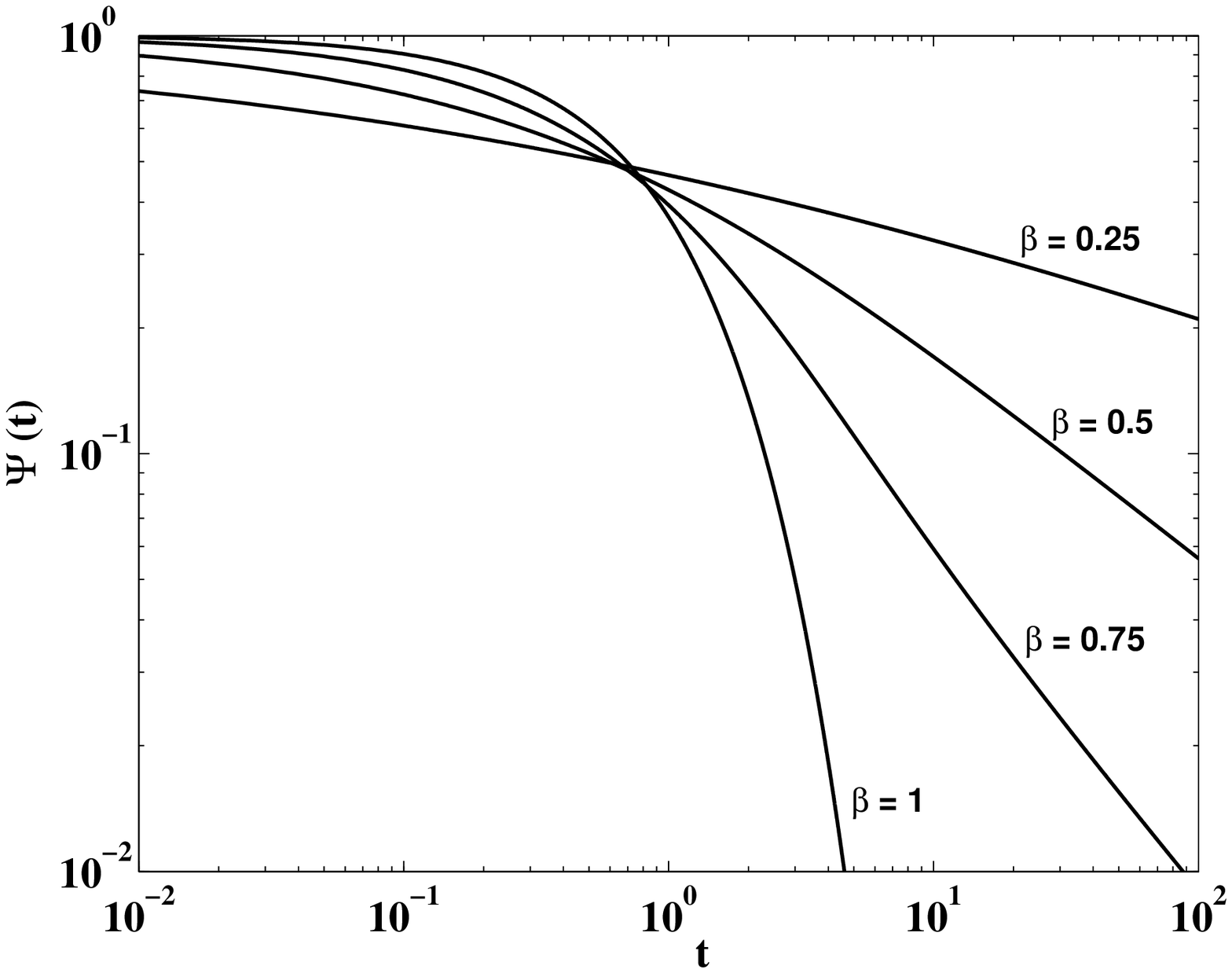}
\includegraphics[width=.52\textwidth]{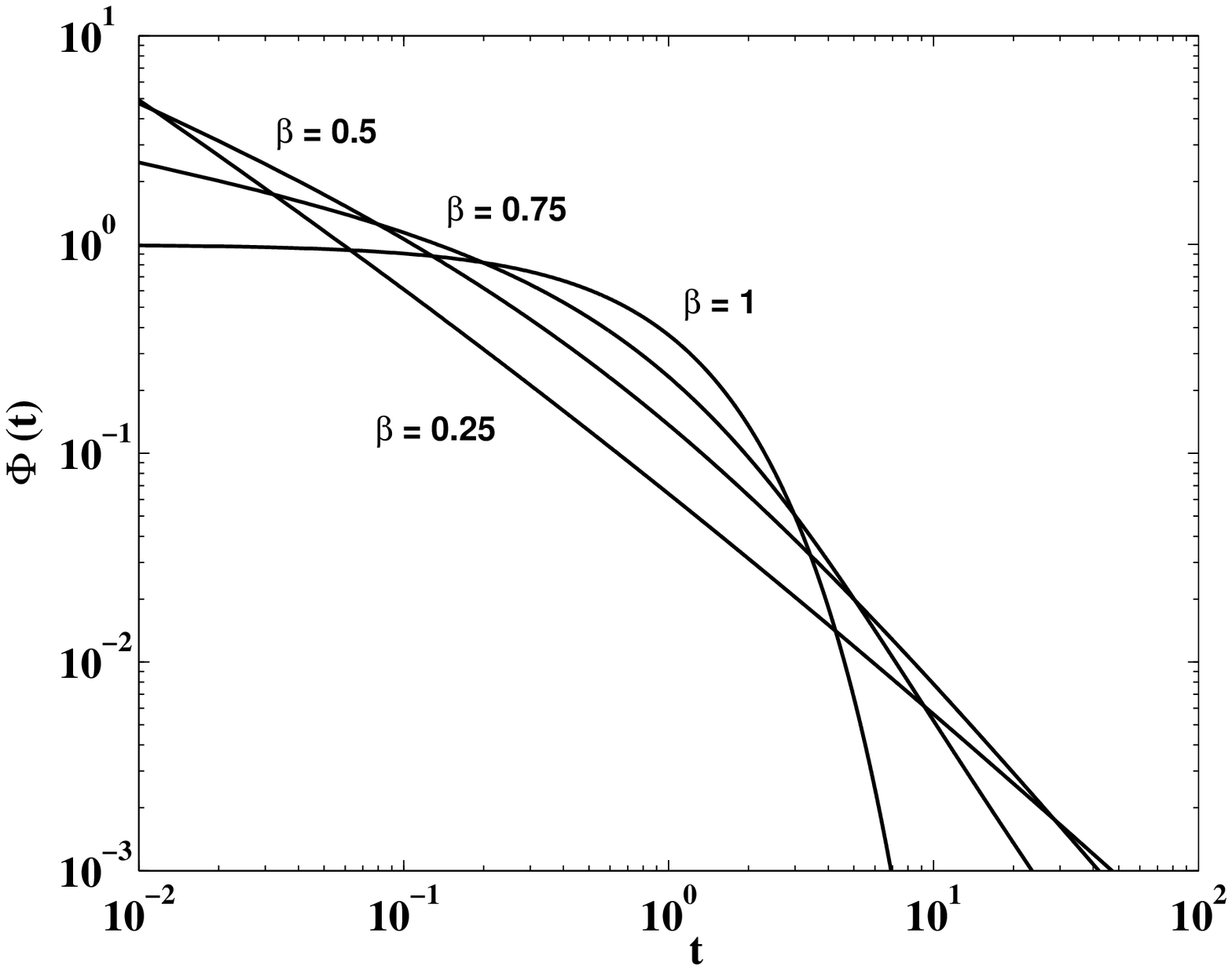}
\caption{The functions $\Psi(t)$ (left) and
$\phi(t)$ (right)
versus $t$ ($10^{-2}<t<10^{2}$) for the renewal  processes
of  Mittag-Leffler  type
with $\beta = 0.25, 0.50, 0.75, 1$.}
\end{figure}
%%%%%%%%%%%   THE END OF THE FIGURE 1 WITH PSFIG %%%%%%%%%%%%%%%%%%%%%%%%%

% \begin{figure}%% [H]
% \includegraphics[width=.48\textwidth]{fmFCAA04_1L.eps}
% \includegraphics[width=.48\textwidth]{fmFCAA04_1R.eps}
% \caption{\small  The functions $\Psi(t)$ (left) and $\phi(t)$ (right)
% versus $t$ ($10^{-2}<t<10^{2}$) for the renewal  process
%  of Mittag-Leffler  type with $\beta = 0.25, 0.50, 0.75, 1$.}
% \end{figure}

%%%%%%%%%%%%%% THE END of FIGURE 1 WITH GRAPHICS  The Mittag-Leffler function

The renewal function of this process
can be deduced from  the Laplace transforms
in (2.19) and  (4.4); we find
$$ \widetilde m(s) = \frac{1}{s^{1+\beta} }\,,
\q \hbox{hence} \q
   m(t) =  \frac{t^\beta}{\Gamma(1+\beta) }
\,, \q t \ge 0\,, \q 0 < \beta \le 1\,. \eqno(4.8)$$
Thus, for $\beta <1$ the renewal function
 turns out super-linear for small $t$ and sub-linear for
large $t$.

For the generalization of Eqs (3.5)-(3.7),   concerning
 the Poisson and the Erlang distributions,
we give the   Laplace transform formula %% pair
$$ \L\{ t^ {\beta \,k}\, E_\beta ^{(k)}
  (-t^\beta ) ;s\} =
        \frac{ k!\, s^{\beta -1}}{(1+s^\beta )^{k+1}}
\,, \q \beta >0 \,, \q k = 0,1, 2, \dots \,.
 \eqno(4.9)$$
with $ {\ds E_\beta ^{(k)}(z) := \frac{d^k}{dz^k}  E_\beta(z)}\,, $
that can be deduced  from the book by Podlubny,
see (1.80) in  \cite{Podlubny BOOK99}.
%% For reader's convenience we report in Appendix B our proof of (4.8)
%% adapted  from \cite{Podlubny BOOK99}.
%%%%%%%%%
Then we get, with $0<\beta <1$ and $k =0, 1, 2, \dots$,
$$ v_k(t) =  \frac{ t^{ k\, \beta}}{k!} \,
  E_\beta^{(k)} (- t^\beta)\,,
\; \hbox{from}    \;
\widetilde v_k(s) =
\widetilde \Psi(s)\, [\widetilde \psi(s)]^k  =
\frac{s^{\beta -1}}{(1+s^\beta)^{k+1}}\,,
  \eqno(4.10)$$
 as  generalization of  the Poisson distribution (with parameter $t$),
what we  call the
{\it $\beta $-fractional Poisson distribution}.
Similarly, with $0<\beta <1$  and $k =1,2, \dots$,
  we obtain  the {\it $\beta$-fractional Erlang} $pdf$ (of order $k \ge 1$):
$$ f_k(t) = \beta \,  \frac{ t^{k\beta-1}}{(k-1)!}
    \, E_\beta^{(k)} (- t^\beta)  \,,
\; \hbox{from}    \;
\widetilde \phi _k(s) =
  [\widetilde \psi(s)]^k  =
\frac{1}{(1+s^\beta)^{k}}\,,
 \eqno(4.11)$$
and the corresponding {\it $\beta$-fractional  Erlang} $pcf$:
$$ F_k(t) = \int _0^t f_k(t')\, dt' =
     1 -  \sum_{n=0}^{k-1}
      \frac{ t^{n \beta}}{n!} \, E_\beta^{(n)} (- t^\beta)  =
 \sum_{n=k}^{\infty}
      \frac{t^{n\beta}}{n!} \, E_\beta^{(n)} (- t^\beta)
\,. \eqno(4.12)$$
%%%%%%%%%%%%%%
% \vfill \eject
%%%%%%%%%%%%%%%%%%%%%%%%%%%%%%%%%%%%%%%%%%
\section*{5. The renewal  process of Wright type}

A possible choice for obtaining an analytically treatable
      variant to the Poisson process
has been suggested by Mainardi et al.
 \cite{Mainardi PhysA00}.
%% see also
%% Barkai \cite{Barkai ChemPhys02} who already treats
%% the corresponding compound process.
%% (where the authors have  paid attention mainly to
%% the application of these processes  to financial markets).
It  is based on the assumption that the {\it waiting-time} $pdf$
$\phi(t)$  is the density of an extremal,  unilateral,
L\'evy stable distribution with
index $\beta\in (0,1)$, which exhibits, as we shall show,
the same power law asymptotics as the corresponding
$pdf$ of  the previous renewal process.
In this case, however, the transition to the limit
$\beta =1$ is singular, and the Poisson process is no longer obtained.
%%%%%%%%
Now  we have for $t \ge 0$,
 $$   \Psi(t) = \cases{
    1-  \Phi_{-\beta,1} \l(- {1\over t^\beta}\r),
      & $\, 0<\beta<1,$ \cr
 \Theta(t) -\Theta(t-1),
     & $\, \ \beta=1,$ \cr}
  \; \hbox{from} \;
  \widetilde \Psi(s) = {1-\e^{\, \ds -s^\beta}\over s},
\eqno(5.1) $$
%%  and
$$ \phi(t) = \cases{ {1\over t}\, \Phi_{-\beta,0} \l(- {1\over t^\beta}\r),
  & $\, 0<\beta<1,$ \cr
 \delta(t-1),
  & $\, \ \beta=1,$ \cr}
 \; \hbox{from} \;
\widetilde\phi(s)  =    \e^{\,\ds -s^\beta},
 \eqno (5.2)$$
where   $\Phi _{-\beta,1}$ and $\Phi _{-\beta,0}$
denote Wright functions\footnote{%%
%% particular cases of the transcendental function
%% (depending on two indices) %% $\Phi_{\lambda,\mu}$,
%% known as the Wright function$\null^{(4)}$.
%%%%%%%%%%%%%%%%%%%%%%%%%  FOOTNOTE  ON WRIGHT FUNCTIONS %%%%%
The  Wright function with parameters $\lambda ,\mu $ is
defined as %%% in the complex plane by the power series
$$ \Phi_{\lambda ,\mu } (z) :=
    \sum_{n=0}^{\infty}\,
   {z^{n}\over n!\Gamma(\lambda \,n+\mu )}\,,
 \q \lambda  >-1\,,\q \mu \in \CC\,, \q z \in
\CC\,.  $$
It is an entire function of order $\rho =1/(1+\lambda)$.
For detailed information on the  functions of Wright type
%% and their Laplace transforms
the reader  may  consult \eg
 \cite{Stankovic 76,GoLuMa 99,KilbasSaigoTrujillo 02,Kiryakova 94,%%
 Mainardi CISM97,Stankovic 70},
and references therein.
We note that in the classical handbook of the Bateman Project
\cite{Erdelyi HTF}, see Vol. 3,  Ch. 18,
%% Miscellaneous Functions, pp. 206-227,
presumably for a misprint, the Wright function is considered
with   $\lambda $  restricted to be non-negative.
}.
%%%%%%%%%%%%%%%%%%%% THE END OF FOOTNOTE  %%%%%%%%%

In view of the presence of the Wright function in (5-1)-(5.2),
%% we agree to refer to this  process as to the
we call this process the {\it renewal process  of Wright type}.

Hereafter,%% like for the Mittag-Leffler-type functions (3.5),
taking  $0<\beta<1$, we quote for $\Psi(t)$ and $\phi(t)$
their  expansions  in powers series of  $t^{-\beta }$
({\it convergent for} $t>0$)
$$  \Psi(t)  =
  {\ds {1\over \pi} }\, {\ds \sum_{n=1}^{\infty}}\,
  (-1)^{n-1} {\ds {\Gamma(\beta\,n ) \over n!}}\,
  {\ds {\sin(\pi \beta n)\over t^{\beta n}}} \,, %% \q t> 0\,,
\q 0<\beta <1\,,
\eqno(5.3)$$
%% and   %%  For the {\it waiting-time} $pdf$  we have
  $$ \phi(t) =
  {\ds {1\over \pi t}}\, {\ds \sum_{n=1}^{\infty}}\,
  (-1)^{n-1} {\ds {\Gamma(\beta\,n+1 ) \over n!}}\,
  {\ds {\sin(\pi \beta n)\over t^{\beta n}}} \,, %%%\q t> 0\,.
\q 0<\beta<1 \,.
\eqno(5.4)$$
We indeed note that the first term of the above series  (5.3)-(5.4)
is identical to the asymptotic representation of the
corresponding functions $\Psi(t)$ and $\phi(t)$ of
the renewal process  of the Mittag-Leffler type,
 see (4.5)-(4.7).
%%%%%%%%%%5
The behaviour near $t=0$ of these functions %%%%$\Psi(t)$ and $\phi(t)$
is provided by the first term  of their  asymptotic expansions
as $t\to 0$, namely  from \cite{Mainardi PhysA00},
 $$
 \Psi(t) \,\sim \,  1- A\, {\ds t^{b/2}} \,
\exp \l(  \,{\ds - B \, t^{-b}}\r),  \q
\phi(t)  \, \sim \, C\, {\ds t^{-c}} \, \exp \l(
  \,{\ds - B \, t^{-b}}\r),
  \eqno(5.5)$$
where\footnote{%%%
%%%%%%%%%%%%%% FOOTNOTE ON MISPRINT IN LUMAPA FCAA01 !! %%%%%%%%
We take this occasion to point out a misprint (sign error) in the paper
\cite{Mainardi LUMAPA01}.
Noting that   the asymptotic representation
as $x \to 0^+$  of the unilateral extremal density of index
$\alpha \in (0,1)$ in  Eq. (4.15) of \cite{Mainardi LUMAPA01}
is equivalent
 to ours of  $\phi(t)$ as $t \to 0$,
 the coefficient $c_1$ (corresponding to our $b$ in (5-5)-(5-6)),
 must be taken with  the {\it minus} sign.}
$$
A =  \l[{1\over 2\pi (1-\beta)\, \beta^{1/(1-\beta )}}\r]^{1/2},\q
   B = (1-\beta )\,\beta ^b,\q
   C = \l[{\beta^{1/(1-\beta )}\over 2\pi (1-\beta)}\r]^{1/2}, \qq
\eqno(5.6)$$
$$
b = {\beta \over 1-\beta }\,, \q
   c = {2-\beta \over 2(1-\beta)}\,.$$

%%%%%%%%%%%% The END of FOOTNOTE
%%%%%%%%%%%%%%%%%%%%%   FIG 2 with PSFIG
\begin{figure}%% [htbp]
 \includegraphics[width=.52\textwidth]{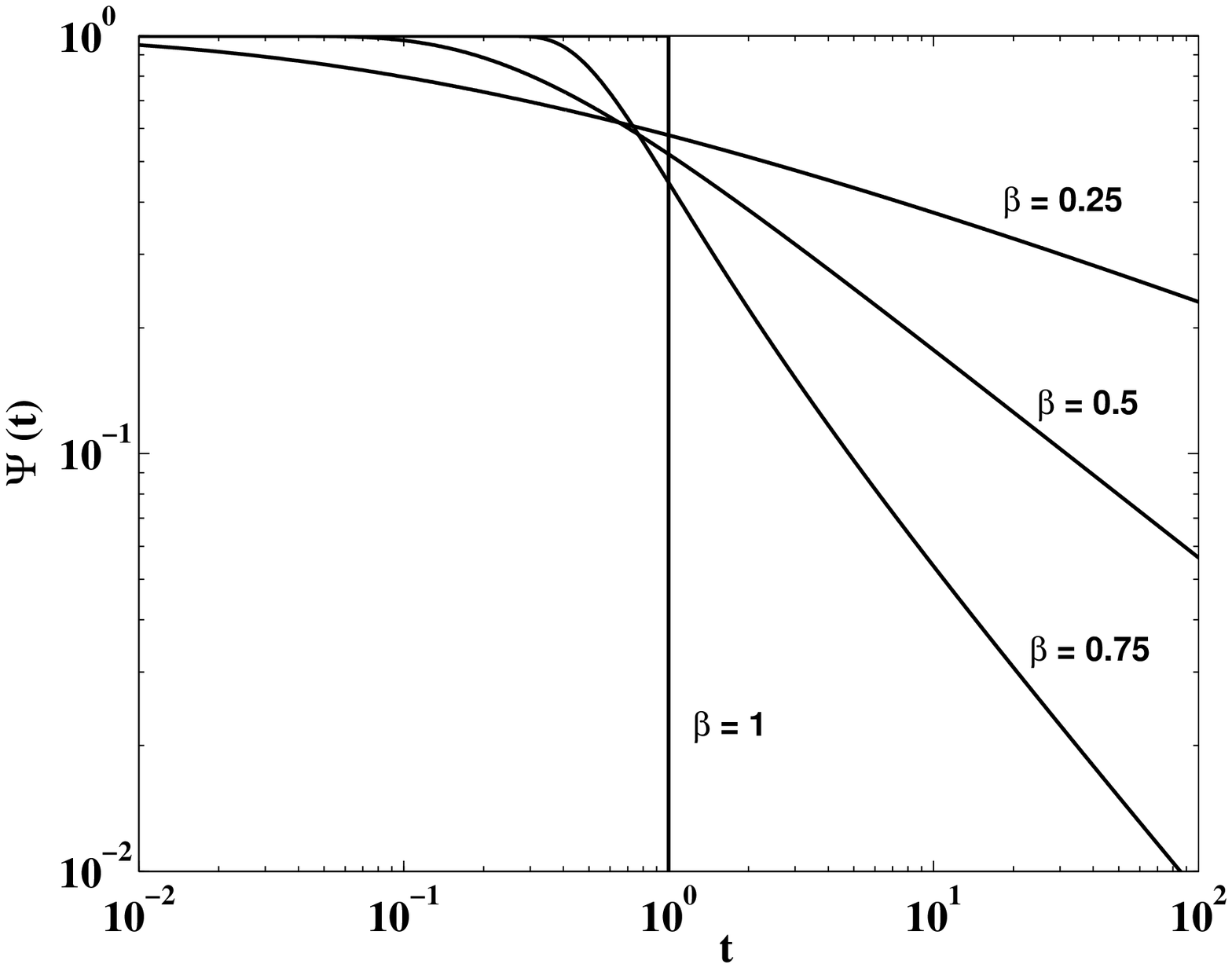}
\includegraphics[width=.52\textwidth]{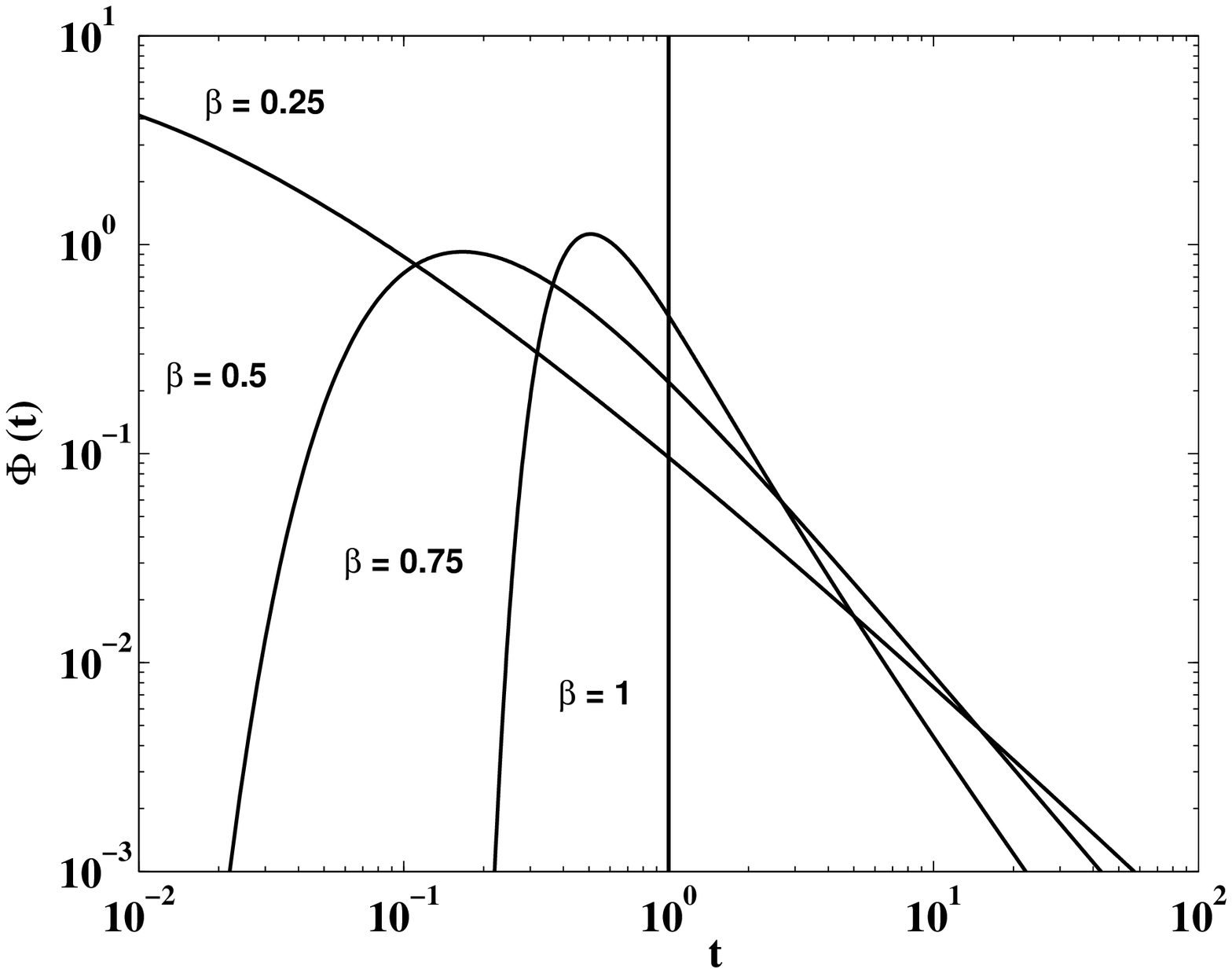}
\caption{ The functions $\Psi(t)$ (left) and
$\phi(t)$ (right)
versus $t$ ($10^{-2}<t<10^{2}$) for the renewal  processes of Wright type
with $\beta = 0.25, 0.50, 0.75, 1$.
For $\beta =1$ the reader would recognize the Box function
(extended up to $t=1$) at left and the delta function
(centred in  $t=1$) at right.}
 \end{figure}
%%%%%%%%%%%   THE END OF THE FIGURE 2 with PSFIG %%%%%%%%%%%%%%%%%%%%%%%%%%%%
% $\begin{figure}%% [H]
% \includegraphics[width=.48\textwidth]{fmFCAA04_2L.eps}
% \includegraphics[width=.48\textwidth]{fmFCAA04_2R.eps}
% \caption{\small  The functions $\Psi(t)$ (left) and $\phi(t)$ (right)
% versus $t$ ($10^{-2}<t<10^{2}$) for the renewal process  of Wright type
% with $\beta = 0.25, 0.50, 0.75, 1$.}
% \end{figure}
%%%%%%%%%%%%%% THE END of FIGURE 2 with GRAPHICS %%%%%%%%
As far as the functions  $v_k(t)$, see (2.14),
%% (that provide the probability that the sum of the first $k$ waiting times
%% is equal to or less than $t$)
are concerned, we have
$$ \widetilde v_0(s) = \widetilde \Psi(s) =
\frac{1 - \e^{\, \ds -s^{\beta}} }{s}\,,\eqno(5.7)$$
so,
 $$  v_0(t)= \Psi(t) = \cases{
    1-  \Phi_{-\beta,1} \l(- {1\over t^\beta}\r),
      & $\, 0<\beta<1,$ \cr
 \Theta(t) -\Theta(t-1),
     & $\, \ \beta=1,$ \cr}
\eqno(5.8) $$
and
$$ \widetilde v_k(s) =  \widetilde \Psi(s)\, [\widetilde \psi(s)]^k
%%   = \frac{1 - \e^{\, \ds s^{\beta}} }{s}\, \e^{\, \ds ks^{\beta}}
 = \frac{ \e^{\, \ds -ks^{\beta}} }{s} - \frac{\e^{\,\ds -(k+1)s^{\beta}}}{s}
\,, \q  k = 1, 2, \dots \,, \eqno (5.9)
$$
from which,
in view of the scaling property of the Laplace transform,
$$ v_k(t) =   \cases{
   {\ds \Phi_{-\beta,1} \l(- {k\over t^\beta}\r)-
\Phi_{-\beta,1} \l(- {k+1\over t^\beta}\r)},
       & $\, 0<\beta<1,$ \cr
 \Theta(t-k) -\Theta(t-k-1),
     & $\, \ \beta=1.$ \cr}
%%%%%\,, \q k = 1, 2, \dots \,.
\eqno(5.10)$$
%%%%%%%%
Let us close this Section  with a discussion on the renewal function
$m(t)$. %%%  for the renewal process of Wright type.
Whereas for the process of Mittag-Leffler type we have an
explicit expression, namely (4.8), we could not find such
one for the process of Wright type in the case
$0 < \beta < 1$ from
the Laplace transforms  (2.19) and (5.2), i.e.
$$  \widetilde m(s) =
\frac{\widetilde \phi (s)}{s\,\l[1-\widetilde \phi (s)\r]}
= \frac{1}{s}\,
   \frac{\e^{\, \ds -s^\beta}}{ 1-\e^{\, \ds -s^\beta}}
  = \frac{1}{s}\, \sum_{k=1}^{\infty}  \e^{\, \ds -k\,s^{\beta}}\,.
\eqno(5.11)$$
In the special case $\beta =1$    we  have
(the term by term inversion is allowed!)
$$ m(t)  = \sum_{k=1}^{\infty} \Theta (t-k) = [t]\,,\eqno(5.12)$$
where
%% $\Theta(t)$   is the Heaviside unit step function
%% with $\Theta(0)= \Theta(0^+) =1$ and
$[t]$ denotes the greatest integer less than or equal to $t$.
%%%%%%%%%

%%%%%%%%%
In the case  $0<\beta <1$  we do not know how to invert  (5.11).
 We can,
however,
see the asymptotics near zero and near infinity and apply  Tauber theory:
%%%%%%%
$ s\to 0$  gives  $\widetilde m(s) = 1/s^{1+\beta }$,
hence,
$$m(t) \sim t^\beta /\Gamma(1+\beta )\q \hbox{for} \q t\to \infty\,,
\eqno(5.13) $$
$ s\to \infty$ gives   $\widetilde m(s) = \exp{(-s^\beta) }/s =
   \widetilde \Phi(s)$,
hence (if this Tauber trick is allowed for such fast decay towards zero),
in view of (5.5)-(5.6),
$$m(t) \sim \Phi(t) \sim
   A\, {\ds t^{b/2}} \,
\exp \l(  \,{\ds - B \, t^{-b}}\r)
\q \hbox{for} \q t\to 0\,. \eqno(5.14)$$

%%%%%%%%%%
\section*{6. The   Mittag-Leffler and Wright processes in comparison}

In this section  we intend to compare  the renewal processes
of Mittag-Leffler and Wright type  introduced in the Sections
4 and 5, respectively.
In this comparison we agree to use  the upper indices $a$ and $b$
to distinguish the relevant functions
characterizing  these processes.
%%%%%%%%%%%

We begin  by pointing out the major differences between
the survival functions
  $\Psi^a (t)$ and $\Psi ^b(t)$, provided by Eqs (4.3) and (5.1), respectively,
for a common index $\beta $ when $0<\beta \le 1$.
These differences, visible from a comparison
of the plots (with logarithmic scales) in the left plates of Figures 1 and 2,
can be easily inferred by analytical arguments,
as  previously pointed out (in a preliminary way)
in the paper by Mainardi et al.
\cite{Mainardi PhysA00}.
%%  where the authors paid attention to
%% the application of these processes  to financial markets.

%% The first difference  concerns  the decreasing behaviour before the
%% onset of the common power law regime:
%% in a Cartesian plot (with linear scales) we can note that,
%% whereas $\Psi^a (t)$
%% starts at $t=0$ vertically (the derivative is $-\infty$)  and is
%% completely monotone,
%% $\Psi^b (t)$    starts at $t=0$ horizontally (the derivative is $0$)
%% and then  exhibits a change from concavity to convexity.
%%%%% in the concavity from downwards to upwards.

Here we stress again the different behaviour
of the two processes in the limit  $\beta \to 1\,$
for $t\ge 0$.
Whereas  $\Psi^a (t)$ and $\phi^a(t)$  tend to the exponential $\exp (-t)$,
$\Psi^b (t) $ tends  to the box function
$\Theta(t) - \Theta(t-1)$
 and  the corresponding {\it waiting-time pdf}
$\phi ^b(t)$ tends to the
shifted Dirac delta  function $\delta(t-1)$.
%% as  directly obtained from the Laplace inversion of Eqs. (5.1)-(5.2)
%% for $\beta=1$. %%  We  fix $\Theta(0)=\Theta(0^+)=1$.
The first of these processes is thus a direct generalization of the Poisson
process  since, for the limiting  value $\beta = 1$, the Poisson process is
recovered. In distinct contrast, the second process changes its character
from stochastic to deterministic: for $\beta = 1$ at every instant $t=n$,
$n$ a natural number, an event happens (and never at other instants)
so that $N(t) = [t]$,
hence trivially $m(t) :=E (N(t)) = [t]$
as we had already found in (5.12) by summation.
%%%%%%%%%%%%%%%%%%%
We refer to this peculiar  counting process as the
{\it clock  process} because  of its similarity with
the tick-tick of a (perfect) clock. We also note that in the limit
$\beta =1$ the densities of the
Mittag-Leffler and  Wright processes have an identical finite
first moment since
$\rho = \int_0^\infty t \, \exp (-t)\, dt =
\int_0^\infty t \, \delta (t-1)\, dt =1$.
%%%%%%%%%%%%%%%%%%%%%%

%%%%%%%%%%%%%%%
It is instructive to consider the special value $\beta =1/2$
because in  cases ($a$) and ($b$)  we have  an explicit
representation of the corresponding survival functions
and waiting-time densities
in terms of well known functions.
%% (the error  and complementary error function).
%%%%%%%%%%

For the renewal process of Mittag-Leffler type with $\beta =1/2$ we have
$$ \Psi^a(t) =  E_{1/2} (-\sqrt{t}) =
    \e^{\ds \, t}\, \hbox{erfc} (\sqrt{t}) =
 \e^{\, \ds t}\, {2\over \sqrt{\pi}}\,
  \int_{\sqrt {t}}^\infty \e^{\, \ds -u^2}\,du  \,,\q t\ge 0\,,
\eqno(6.1)$$
where $ \, \hbox{erfc}\,$ denotes the
{\it complementary error} function, see \eg \cite{AS 65},
% \footnote{%%
%%%%%%%%%%%%%%%%%%%%%%%%%%%%%%%%%
% We remind for $z \in \CC$,
%% $\erf(z)$ and $\erfc(z)$,
% see \eg \cite{AS 65},
% $$ \erf (z) := {2\over \sqrt{\pi}}\,
%  \int_{0}^{z} \e^{\, \ds -\zeta^2}\,d\zeta =
%\sum_{n=0}^\infty  {2^{\ds n}\over (2n+1)!! }\, z^{\ds 2n+1}, \;
%   \erfc (z) := 1- \erf (z),$$
 %% = {2\over \sqrt{\pi}}\, \int_z^{\infty} {\e}^{\ds-\zeta^ 2}\, d\zeta
%% It is useful to  also recall the asymptotic expansion of $\erfc(z)$:
%$$ \erfc (z) \sim
% \rec{\sqrt{\pi}}\,  {\;\;{\rm e}^{\ds -z^2} \over z}\,
%\l(1 + \sum_{m=1}^\infty  {(-1)^{\ds m}\over m!}
%\, {(2m)!\over  (2z)^{\ds 2m}}\r)\,,
%   \;  |z| \to \infty\,, \; |{\rm arg}\, z| < {3\pi\over 4}\,.$$
%},
%%%%%%%%%%%%%%%%%%% END OF FOOTNOTE on ERROR FUNCTIONS
 and
$$    \phi^a(t) = - {d \over dt} E_{1/2} (-\sqrt{t}) =
  \frac{1}{ \sqrt{\pi t}} -  \e^{\ds \, t}\, \hbox{erfc} (\sqrt{t})\,,
\q t \ge 0\,.\eqno(6.2)$$

%%%%%%%%
For the renewal process of Wright type with
 $\beta ={1/ 2}$
 we obtain for $t \ge 0$, %% the explicit expressions
$$   \Psi^b(t) = 1- \Phi_{-1/2, 1} (-1/\sqrt{t}) =
 1 - {\ds \hbox{\erfc} \l({1\over 2\sqrt{t}}\r)} =
{\ds \hbox{\erf} \l({1\over 2\sqrt{t}}\r)}    \,,
 \eqno(6.3)$$
 and
$$ \phi^b(t) = {1\over t}\,\Phi_{-1/2, 0} (-1/\sqrt{t}) =
{1\over 2\sqrt{\pi}}\, t^{-3/2}\,\exp \l(- {1\over 4t}\r)\,.
\eqno(6.4)  $$
We observe that for this particular value of $\beta$
 the  expression for the density\footnote{%%} 
%%%%%%%%%%%
We point out %% that
$
{\ds \L\l \{\phi^b(t) =
{t^{-3/2}\over 2\sqrt{\pi}}\,\exp \l(- {1\over 4t}\r); s \r \}
= \exp (-s^{1/2})}.$
This Laplace transform pair
was noted by L\'evy with respect to the unilateral
stable density of order $1/2$  and later, independently, by Smirnov.
In the probability literature the distribution corresponding to
this $pdf$ is  known   as  the
{\it L\'evy-Smirnov}  stable distribution.
We take this occasion to point pout a misprint in the paper
\cite{Mainardi PhysA00}. There, in Eq. (3.25) giving the
expression of the L\'evy-Smirnov density,
 the factor $2$ in front of $\sqrt \pi$ was missed.}
%%%%%%%%%%%%%%%%%%%%%% 
 is obtained not only
 from the the sum of the convergent series (5.4)
but also exactly from its  asymptotic representation
for $t \to 0$, see  (5.5)-(5.6).

We easily note the common asymptotic (power-law) behaviour
of the survival and density functions as $t \to \infty$
in the cases ($a$) and ($b$), that is,
indicating by the index $\infty$ that we mean the
leading asymptotic term:
%  $$     \Psi_\infty(t) = \frac{t^{-1/2}}{\sqrt{\pi}}\,, \q
%  \phi_\infty (t) = \frac{t^{-3/2}}{2 \sqrt{\pi}}\,.\eqno(6.5)$$
 $$     \Psi_\infty(t) = \frac{t^{-1/2}}{\sqrt{\pi}}\,, \eqno(6.5)$$
$$ \phi_\infty (t) = \frac{t^{-3/2}}{2 \sqrt{\pi}}\,.\eqno(6.6)$$
%%%%%%%%%%
It is now interesting to   compare numerically
the  survival functions
and the  density functions
of the two processes,
that is (6.1)-(6.2) %%  with $\beta =1/2$ and %% (2.43)
with (6.3)-(6.4), respectively.
%%%%%%%%%%%%%%%%%%%%%%%%%%%%%% %%%%%%%%%%%%%%%
We find it worthwhile to  add in the comparison
their asymptotic expression (6.5)-(6.6)
and also  the  corresponding  functions of the
Poisson process, namely (3.1)-(3.2).

In Figure 3, we display the plots concerning the above comparison
for the survival functions (left plate)
and for the density functions (right plate)
adopting the continuous line for the exact functions
of the two long-memory processes and  the dashed line
both for the  asymptotic power-law functions
and for the exponential functions of the  Poisson process.
%%%%%%%%%
The comparison is made
by using logarithmic scales for $10^{-1} \le t \le 10^1$,
just before the
onset of the common power law regime;
 we note that,
at least  for the case $\beta =1/2$ under consideration,
the Wright process reaches the asymptotic power-law regime
a little bit earlier than the corresponding Mittag-Leffler process.
%%%%%%%%%%%% FIGURES %%%%%%%%%%
\begin{figure}%% [htbp]
 \includegraphics[width=.52\textwidth]{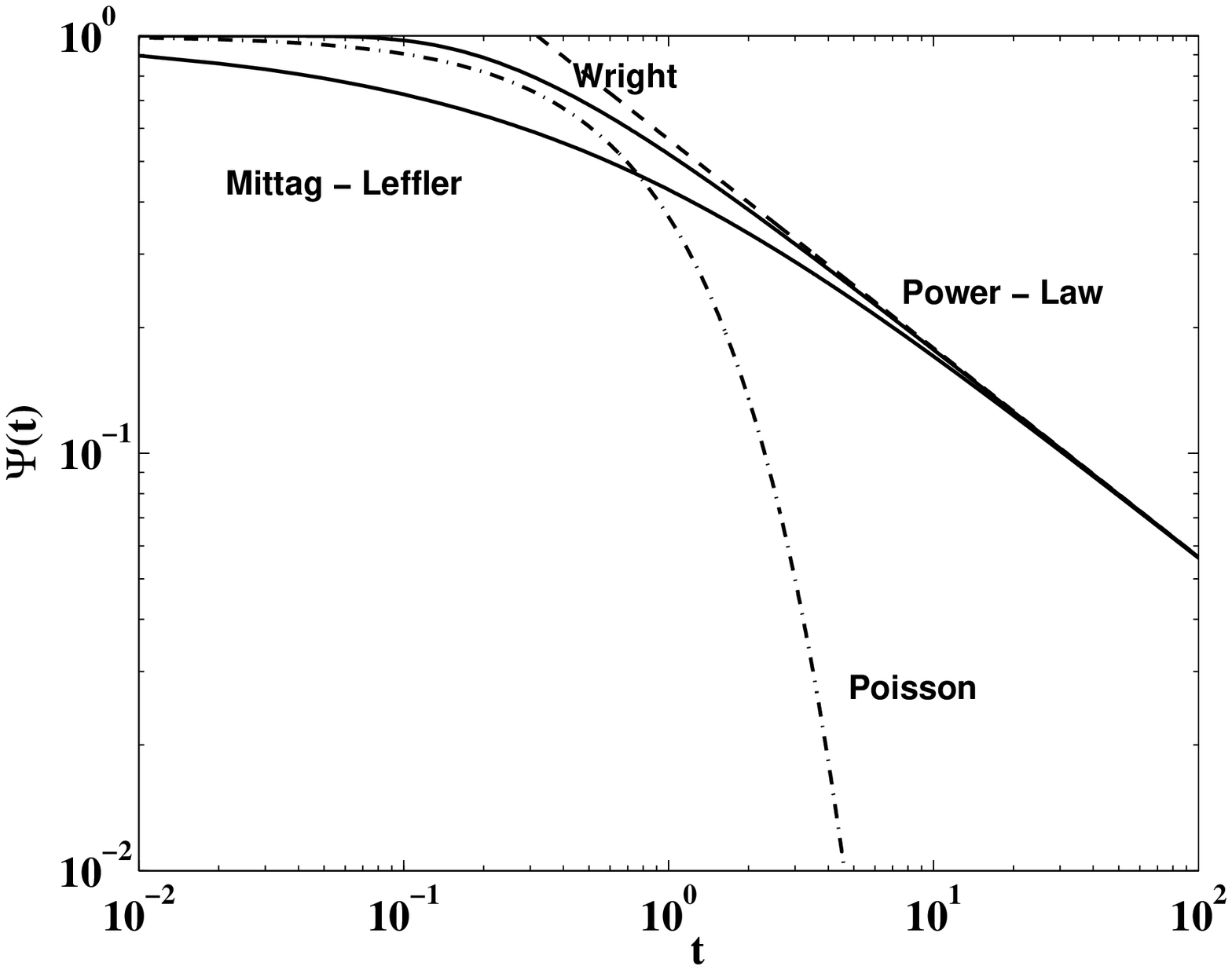}
\includegraphics[width=.52\textwidth]{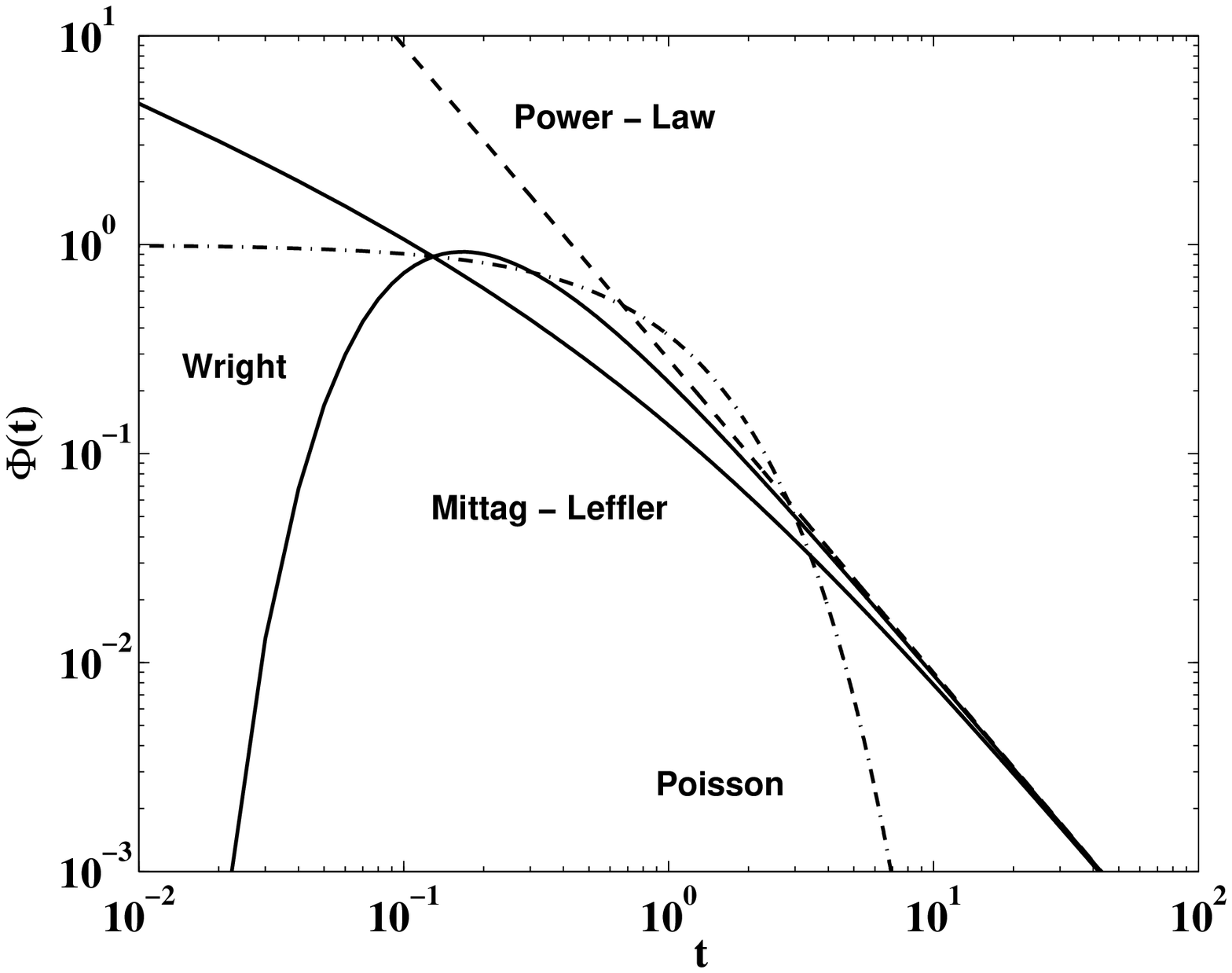}
\caption{Comparison  versus time of the survival functions
 $\Psi(t)$ (left plate) and of the corresponding probability
  densities $\phi(t)$ (right plate).}
%% for the  processes of Mittag-Leffler and Wright type with $\beta=1/2$.
\end{figure}

%%%%%%%%%%%
% \begin{figure}%% [H]
% \includegraphics[width=.48\textwidth]{fmFCAA04_3L.eps}
% \includegraphics[width=.48\textwidth]{fm_FCAA0_3R.eps}
%  \caption{\small Comparison  versus time of the survival functions
% $\Psi(t)$ (left plate) and of the corresponding probability
%  densities $\phi(t)$ (right plate)}
%% Comparison  versus time of the survival functions
%% $\Psi(t)$ (left plate) and of the corresponding probability
%% densities $\phi(t)$ (right plate)
%% of the  Mittag-Leffler and Wright  processes (in continuous line) with
%% their common asymptotic power-laws and the
%% corresponding functions of the Poisson process (in dashed line),
%% for $\beta =1/2$.}
% \end{figure}
%%%%%%%%%%   END OF FIGURE 3 WITH GRAPHICS
In Tables I and II we show %% in 5 columns<
for some values of time $t$
the corresponding values of the above  survival and density functions
in comparison,
abbreviating  power law by P-L and  Mittag-Leffler by M-L.
%%%%%%%%%\vfill\eject
%%%%%%%%%%%% START OF THE TABLE I Survival Functions %%%%%%%%
% $\null$
% \vskip -1.0truecm
\begin{center}
% \vskip -1.0truecm
\begin{tabular}{|c||c||c|c||c|}  %% 5 columns
\hline
{Time}
& P-L (6.5)& M-L (6.1) & Wright (6.3) & Poisson (3.1) \\
 \hline
0.1 & $1.78\,10^{0} $ & $7.24\,10^{-1}$ & $9.74\,10^{-1}$ & $9.05\,10^{-1}$ \\
0.5 & $7.98\,10^{-1}$ & $5.23\,10^{-1}$ & $6.83\,10^{-1}$ & $6.07\,10^{-1}$ \\
1   & $5.64\,10^{-1}$ & $4.28\,10^{-1}$ & $5.21\,10^{-1}$ & $3.68\,10^{-1}$ \\
2   & $3.99\,10^{-1}$ & $3.36\,10^{-1}$ & $3.83\,10^{-1}$ & $1.35\,10^{-1}$ \\
5   & $2.52\,10^{-1}$ & $2.32\,10^{-1}$ & $2.48\,10^{-1}$ & $6.74\,10^{-3}$ \\
10  & $1.78\,10^{-1}$ & $1.71\,10^{-1}$ & $1.77\,10^{-1}$ & $4.54\,10^{-5}$ \\
20  & $1.26\,10^{-1}$ & $1.23\,10^{-1}$ & $1.26\,10^{-1}$ & $2.06\,10^{-9}$ \\
50  & $7.98\,10^{-2}$ & $7.90\,10^{-2}$ & $7.97\,10^{-2}$ & $-$\\%%$1.93\,10^-22$  \\
100 & $5.64\,10^{-2}$ & $5.61\,10^{-2}$ & $5.64\,10^{-2}$ & $-$\\%%$3.72\,10^-44$
\hline
\end{tabular}
\vskip 0.20truecm
{\bf Table I}: Comparison among
the survival functions at different times.
\end{center}

%%%%%%%%%%%%%  END OF THE TABLE I %%%%%%%%%

%%%%%%%%%%%% START OF THE TABLE II  Density functions  %%%%%%%%

\begin{center}
%  \vskip -0.05truecm
\begin{tabular}{|c||c||c|c||c|}  %% 5 columns
\hline
{Time}
& P-L  (6.6) & M-L (6.2) & Wright (6.4) & Poisson (3.2)  \\
 \hline
0.1 & $8.92\,10^{0} $ & $1.06\,10^{-1}$ & $7.32\,10^{-1}$ & $9.05\,10^{-1}$ \\
0.5 & $7.98\,10^{-1}$ & $2.75\,10^{-1}$ & $4.84\,10^{-1}$ & $6.07\,10^{-1}$ \\
1   & $2.82\,10^{-1}$ & $1.37\,10^{-1}$ & $2.20\,10^{-1}$ & $3.68\,10^{-1}$ \\
2   & $9.97\,10^{-2}$ & $6.27\,10^{-2}$ & $8.80\,10^{-2}$ & $1.35\,10^{-1}$ \\
5   & $2.52\,10^{-2}$ & $2.00\,10^{-2}$ & $2.40\,10^{-2}$ & $6.74\,10^{-3}$ \\
10  & $8.92\,10^{-3}$ & $7.83\,10^{-3}$ & $8.70\,10^{-3}$ & $4.54\,10^{-5}$ \\
20  & $3.15\,10^{-3}$ & $2.94\,10^{-3}$ & $3.11\,10^{-3}$ & $2.06\,10^{-9}$ \\
50  & $7.98\,10^{-4}$ & $7.75\,10^{-4}$ & $7.94\,10^{-4}$ & $-$\\%%$1.93\,10^-22$  \\
100 & $2.82\,10^{-4}$ & $2.78\,10^{-4}$ & $2.81\,10^{-2}$ & $-$ \\ %%$3.72\,10^-44$
\hline
\end{tabular}
 \vskip 0.20truecm
{\bf Table II}: Comparison among
the density functions at different times.
\end{center}
%%%%%%%%

We finally consider the functions $v_k(t)$
for the two long-memory processes.
%%%%%%%%
For the functions $v_k^a(t)$ of the Mittag-Leffler process,
see (4.10),
we can take profit of the recurrence relations for repeated integrals
of the error functions, see   \eg \cite{AS 65},  \S 7.2, pp 299-300,
to compute the derivatives of the Mittag-Leffler functions
in Eqs. (4.10).
%% For this purpose we recall:
We recall for $n= 0, 1, 2, \dots\,,$
$$ \frac{d^n}{dz^n} \l( \e^{z^2}\, \erfc (z) \r) =
  (-1)^n \, 2^n\, n!\, \e^{z^2}\, I^n \, \erfc (z)\,, \eqno(6.7)$$
where
$ I^n \, \erfc (z) =\int_z^\infty I^{n-1} \,\erfc (\zeta)\, d\zeta$
and
$I^{-1} \, \erfc (z) =  2 \, \exp (-z^2) /\sqrt{\pi}$.
%%%%%%%

For the  Wright process, see (5.10),
we have
$ v^b_0(t) = \Psi^b(t)$ as in (6.3)
%%  the functions $v^b_k(t)$ read %%% turn out to be
%% $$ v_0(t) =  {\ds \hbox{\erf} \l({1\over 2\sqrt{t}}\r)}\,, $$
and
$$ v^b_k (t) =
\l[  {\ds \hbox{\erfc} \l({k\over 2\sqrt{t}}\r)}-
{\ds \hbox{\erfc} \l({k+1\over 2\sqrt{t}}\r)}\r]\,, \q k=1,2, \dots\,.
\eqno(6.8)$$
%%%%%%%%%%%%%\vfill\eject
%%%%%%%%
\section*{7. Renewal processes with reward}

The renewal process can be accompanied  by reward that means
that at every renewal instant a space-like variable makes a  random jump
from its actual position to a new point in "space".
"Space" is here  meant in a very general sense.
In the insurance business, e.g., the renewal points
are instants where the company receives a payment or must
give away money to some claim of a customer,
so space is money.
In such process occurring in time and in space,
also referred to as {\it compound renewal process},
the probability distribution of jump widths is as relevant
as that of waiting times.

Let us denote by $X_n$  the jumps occurring at instants $t_n\,,$
$\, n = 1,2,3,\dots,$ and
assume that they  are $iid$ (real, not necessarily positive)
 random variables
with  probability density $ w(x)$,  independent of  the
{\it waiting time} density $\phi (t)$.
In a physical context the $X_n$s represent
the jumps of a diffusing particle (the walker),
and the resulting random walk model is known
as {\it continuous time random walk}, abbreviated as CTRW\footnote{%%
%%%%%%%%%%%%%%%%%%%   FOOTNOTE ON CTRW %%%%%%%%
The name CTRW became
popular in physics after
 Montroll, Weiss and Scher (just to cite  the pioneers)
in the 1960's and 1970's published a celebrated series
of papers on random walks for modelling
diffusion processes on lattices, see \eg
\cite{MontrollScher 73,MontrollWeiss 65}, and
the book by Weiss \cite{Weiss BOOK94} with  references therein.
CTRWs are rather good and general phenomenological models for diffusion,
including processes of anomalous transport.
However,  it should be noted that
the idea of combining a stochastic process for waiting times between
two consecutive events and another stochastic process which associates
a reward or a claim to each event dates back at least to the first half
 of the twentieth century with the so-called
 Cram\'er--Lundberg model for insurance risk, see
\cite{Embrechts 01} for a review.},
%%%%%%% THE END OF FOOTNOTE on CTRW  %%%%%%%%%%%%%%%%%%%%
in that the  waiting time is assumed to be
a random variable
with a {\it continuous} probability distribution function.

The position $x$ of the walker at time $t$ is
$$
x(t) = x(0) + \sum_{k=1}^{N(t)} X_k\,,
\eqno(7.1)$$
with  $N(t)$ as in (2.5).
%%% \sum_{k=1}^{n: t_n \leq t} X_k
%% The upper limit has been replaced because of
%% the definition of $N(t)$ !!!!!!!!!!!!!!!
%%%%%%%%%%%

Let us now denote by $p(x,t)$ the probability density
(density with respect to $x$)  of finding the random walker
at the position $x$ at time instant $t$.  We assume
the initial condition $p(x,0) = \delta (x)\,, $  meaning that the
walker is initially at the origin,  $x(0)=0\,. $
We look for the evolution equation for $p(x,t) $
%%% that we shall call   the {\it master equation}
 of the compound renewal process.

Based upon the previous probabilistic arguments
 %%%  the required {\it master equation} reads
we arrive at
$$   p(x,t) =  \delta (x)\, \Psi(t) +
   \int_0^t   \phi (t-t') \, \l[
 \int_{-\infty}^{+\infty}  w(x-x')\, p(x',t')\, dx'\r]\,dt'
 \,, \eqno(7.2) $$
called the {\it integral equation of the CTRW}.
%%%%%%%%%%% \vfill\eject%%%%%%%%%%%%%
%% In fact,
From Eq. (7.2) we recognize the role of the
{survival probability}
$\Psi(t)$ and of the densities $\phi (t)\,,\,  w(x)\,.$
The first term in the RHS of (7.2) expresses the persistence
(whose strength decreases with   increasing time)
of the initial position $x=0$.
The second term (a space-time convolution) gives
the contribution to $p(x,t)$ from the walker sitting
in point $x' \in \RR$ at instant $t' < t$ jumping to
point $x$ just at instant $t\,,$  after stopping (or waiting) time
$t-t'\,. $

The integral equation (7.2) can be solved
by using the machinery of the transforms of Laplace
and Fourier.
Having introduced the notation for the Laplace transform
in sec. 1, we now quote
our notation for the Fourier transform,
%% $-\infty <\kappa <\infty $  are
%% $$ \widetilde{g}(s)= \int_{0}^{\infty} \e^{-st}\, g(t)\, dt
$ \F\{f(x); \kappa\} = \widehat{f}(\kappa)=\int_{-\infty}^{+\infty}
\e^{i \kappa x} \,f(x) \,dx $ ($\kappa \in \RR$),
 and for the corresponding Fourier  convolution of
(generalized) functions
%% $$(g_{1}(t)*g_{2}(t)= \int_{0}^{\infty}
%%  g_{1}(t')\,g_{2}(t-t')\, dt' \q $$
$$ \l(f_{1}\,*\, f_{2}\r)(x)= \int_{-\infty}^{+\infty}
f_{1}(x')\,f_{2}(x-x')\, dx'\,.$$
%% $$
Then, applying the transforms of Fourier and Laplace in succession to
(7.2) and using the well-known
operational rules, we
arrive at %%% the relation
%% $$
%% \widehat{\widetilde{p}}(\kappa,s)=\frac{1-\widetilde{\phi}(s)}{s}\, +
%% \widetilde{\phi}(s)\,\widehat{w}(\kappa)
%% \widehat{\widetilde{p}}(\kappa,s)  \, .
%%\eqno(7.3)$$
%% which   leads to
the famous Montroll-Weiss equation, see \cite{MontrollWeiss 65},
$$
 \widehat{\widetilde{p}} (\kappa,s) = %%% \frac{1-\tilde{\phi}(s)}{s} \,
\frac{\widetilde\Psi(s)}{1- \tilde{\phi}(s)\,\widehat{w}(\kappa)}
\,.
\eqno(7.3)
$$
As pointed out in \cite{GAR Vietnam03},
this equation can alternatively be derived from the Cox formula, see
    \cite{Cox RENEWAL67}  chapter 8 formula (4),
describing the process as    subordination of a random
    walk to a renewal process.
%%%%%%%%%\vfill\eject%%%%%%%%%%%%%

By inverting the transforms one can, in
    principle, find the evolution $p(x,t)$  of the sojourn density for
    time $t$    running from zero to infinity.
In fact, recalling that $|\widehat w(\kappa)| < 1$ and
$|\widetilde\phi (s)| < 1$,
if $\kappa \not= 0$ and
$s \not= 0$, Eq. (7.3) becomes
$$
\widetilde{\widehat p}(\kappa, s) =
\widetilde \Psi(s)\, \sum_{k=0}^{\infty}
[\widetilde \phi (s) \, \widehat w(\kappa)]^k =
 \sum_{k=0}^{\infty} \widetilde v_k(s)\, \widehat w_k (\kappa)
 \,,
\eqno(7.4) $$
and %% inverting the  transforms
we promptly obtain the series representation
%% of the  {\it fundamental solution}  of the CTRW
$$
p(x,t) = \sum_{k=0}^{\infty} v_k(t)\, w_k (x)
=    \Psi(t)  \, \delta (x) + \sum_{k=1}^{\infty} v_k(t)\, w_k (x)\,,
\eqno(7.5)$$
where the functions  $v_k(t)$
and $w_k(x)$ are obtained by repeated convolutions
in time  and in space,
$ v_k(t) = (\Psi * \phi^{*k})(t)$, see also (2.14),
%%   \phi_k(t) = \l(\phi^{*k}\r)(t)\,;\qq
and $\, w_k(x) = (w ^{*k})(x)$, respectively.
In particular,
$ v_0(t) = (\Psi*\delta) (t)= \Psi(t)$, $\, v_1(t) = (\Psi * \phi)(t)$,
$\, w_0(x) = \delta(x), \; w_1(x) = w(x).$
%% Details %%  on the above approach to CTRW
%% are found in  papers of our group,
%% see \eg  \cite{Mainardi_Bonn00,GorMai_INDIA03,Scalas_PRE04},
%% and references therein.
%%%%%%%%
In the R.H.S  of Eq (7.5) we  have isolated
the first singular term  related to the initial condition
$p(x,0) =  \Psi(0) \, \delta (x) =\delta (x)$.
%%%%%%
%% \vfill\eject
%%%%%%%%

A special case of the integral  equation (7.2) is obtained for
the {\it compound Poisson process} where $\phi (t) = \e^ {-t}$
(as in (3.2), with $\lambda =1$ for simplicity).
Then, the  corresponding equation reduces after some
manipulations, that best are carried out in the Laplace-Fourier domain,
  to the {\it Kolmogorov-Feller equation}:
$$
  \frac{\d }{\d  t}\,p(x,t)= -p(x,t)+\int_{-\infty}^{+\infty}
w(x-x')\, p(x',t)\, dx' \, ,
\eqno(7.6)
$$
which is the {\it master equation of the compound Poisson process}.
In view of Eqs (3.5) and (7.5) the solution reads
$$p(x,t) = \sum_{k=0}^{\infty} \frac{ t^k}{k!} \,
{\e}^{- t} \,w_k (x)\,.
\eqno(7.7)
$$
When the survival probability is the Mittag-Leffler function
introduced in    (4.3),  the master equation
of the corresponding {\it compound process of Mittag-Leffler type}
can be shown to be
   $$
  \null_tD_*^\beta \,p(x,t)= -p(x,t)+\int_{-\infty}^{+\infty}
w(x-x')\, p(x',t)\, dx' \,, \q 0<\beta <1\,,
\eqno(7.8)
$$
where  $\,_tD_*^\beta$ denotes the time fractional
derivative   of order $\beta $  in the Caputo sense.
%%% briefly introduced in the Appendix.
For a (detailed) derivation of Eq (7.8 ) we refer  to the paper
by Mainardi et al. \cite{Mainardi PhysA00}, in which the results
have been obtained by an  approach different
from that adopted in a previous paper
by Hilfer and Anton \cite{HilferAnton 95}.
%%%%%%%%%%
In this case, in view of Eqs (4.10) and (7.5), the solution
of the {\it  master equation} (7.8) reads:
$$ p(x,t) =   E_\beta (- t^\beta)  \,\delta  (x) +
\sum_{k=1}^{\infty} \frac{ t^{\beta k}}{k!} \,
  E_\beta^{(k)} (- t^\beta)  \,w_k (x)\,,
\q 0<\beta <1\,.
\eqno(7.9)
$$
%%%%%%%%
When the survival probability is the Wright function
introduced in    (5.1),   we resist the temptation to work out an
explicit form for the the master equation
of the corresponding {\it compound process of Wright type},
but  we content ourselves
to give its series solution derived from  (7.5)
using Eqs (5.8),(5.10),
as
 $$ \qq \qq \qq \qq \qq
p(x,t) = \l[1-\Phi_{-\beta, 1}\l(-\frac{1}{t^\beta }\r)\r]\,\delta(x) \; +
\qq \qq \qq  \qq \eqno(7.10)$$
$$ \sum_{k=1}^{\infty}
\l[\Phi_{-\beta, 1}\l(-\frac{k}{t^\beta }\r) -
\Phi_{-\beta,1}\l(-\frac{k+1}{t^\beta }\r)\r]\, w_k(x)\,, \q 0<\beta <1
 \,. $$
Recalling what we have observed in Section 6 on the contrasting
behaviours of the Mittag-Leffler process and the Wright process in
the limit $\beta = 1$ we recognize that {\it with rewards} the first
of these processes becomes the {\it compound Poisson process}
((7.9) with $\beta=1$ reduces to (7.7)), whereas
the second goes over into the  process,
whose sojourn $pdf$ is directly obtained from the Laplace inversion of
Eqs (5.7) and (5.9)  for $\beta =1$.
Recalling the natural convention $w_0(x)=\delta (x)$,
we obtain
$$ p(x,t) =
 \sum_{k=0}^{\infty} \l[\Theta (t-k)- \Theta(t-k-1)\r]\,w_k(x)\,.\eqno(7.11)$$
This simply means:
$p(x,t) = w_k(x)$ for
     $k \le t < k+1$  ($k =0,1,2, \dots $),
or, more concisely, $ p(x,t) = w_{[t]}(x)$.
%%%%%%%%%%%%%%%%%
% $$ p(x,t) = w_k(x) \q \hbox{for} \q
%%     k \le t < k+1, \q k =0,1,2, \dots\,, \eqno(7.12)$$
%% or, more concisely, $ p(x,t) = w_{[t]}(x)$.
This limiting process is nothing but  a  simple
{\it random walk, discrete  in time and Markovian on the set of integer
time-values}:
  at every instant $[t]$ a jump   occurs with $pdf$ $w(x)$, hence
   $[t]$ jumps have occurred exact up to instant $t$, which  leads to
   $[t]$-fold convolution $w_{[t]}(x)$.
 We have rigorously investigated such random walks in
 \cite{GorMai CHEMNITZ01}.
%%  R. Gorenflo  and  F. Mainardi,
%%  Random walk models approximating symmetric
%% space fractional diffusion  processes,
%%  in  J. Elschner, I.  Gohberg  and  B. Silbermann (Editors),
%% {\it Problems} {\it in Mathematical} {\it Physics},
%%   Birkh\"auser Verlag, Basel, 2001, pp. 120-145.
%% [Series "Operator Theory: Advances and Applications" No 121]
%%% Chemnitz.tex paper %% Siegfried Pr\"ossdorf Memorial Volume,
%% Proceedings of the 11-th TMP Conference
%% Let us cite as an example for discrete in time but continuous in
%% space the {\it Chechkin-Gonchar random walk}.  "
However, in the framework of the CTRW theory  (where
the resulting  random walk, discontinuous or continuous in space,
is implicitly intended to be {\it continuous in time}),
this case appears a peculiar
(non-Markovian) process, as discussed  by  Weiss \cite{Weiss BOOK94}
in his book  at p. 47, see also \cite{Mainardi  PhysA00}.
%% In a different framework appears as  a peculiar example  of a
%% case has been considered
%%%%  as an example of a non-Markovian process.
%% Weiss writes:
%% "An illustrative example in which memory is important is furnished
%% by the case $\phi(t) = \delta(t-1)$, where, at any time $t$, the time
%% to the immediately following time depends only on the time of the
%% preceding step, and not on the current time."
%%   On the other hand, WEISS says:
%% ".... the exponential form of $psi(t)$ is "memoryless".
%% This property is not shared by any other form of $\psi(t)$.
%%  Since the probability of the time of a step
%% is INDEPENDENT of any other variable, the process is NON-MARKOVIAN"
%%%%%%%%%%%%%%%
% \vfill\eject
%%%%%%%%%%%%%
% \sect{8. Time fractional diffusion equation as limit of the CTRW equation}

\section*{8. Time fractional diffusion  as limit of  CTRW}

There are several ways to generalize the classical diffusion equation
by introducing space and/or time derivatives of fractional order.
We mention here   the seminal papers
by Schneider \& Wyss \cite{SchneiderWyss 89} of 1989
for the time fractional diffusion equation and
the likewise influential paper  by Saichev \& Zaslavsky 
\cite{Saichev 97} of 1997 for diffusion fractional in time as well in space.
In the recent literature  several authors
stress the viewpoint
of {\it subordination}, see 
\eg \cite{Baeumer-Meerschaert FCAA01,%%
Bazhlekova FCAA00,Mainardi FCAA03,M3 PRE02sub,Wyss-Wyss FCAA01},
and special attention 
is being paid to diffusion equations
with {\it distributed  orders} of fractional
temporal or/and spatial differentiation,
see \eg
\cite{Caputo FCAA01,Chechkin FCAA03,GorMai CARRY04,Sokolov 04}.
%%%%%%%%%%%
 The transition from CTRW to such generalized types of  diffusion 
has been  investigated by different methods,
not only from mathematical interest but also  due to applications
 in Physics, Chemistry
(for  generalized Fokker-Planck and Liouville equations
 see \eg
\cite{Barkai ChemPhys02,Barkai PRE00,Hilfer PhysA03,HilferAnton 95,%%%
%% Metzler PRE98,
Sokolov PhysicaA01}),
and other Applied Sciences including  Economics and Finance
(for  {\it financial markets} see
\eg \cite{Gorenflo 01,Mainardi PhysA00,Scalas PhysA00}).
%% but without devoting particular attention to the scaling
%% necessary, in our opinion, to   get in ra rigorous way
%% the required transition.
For a good list of references on these topics we recommend
the review papers by Metzler and Klafter
\cite{MetzlerKlafter PhysRep00,MetzlerKlafter JPhysA04}.

Gorenflo and Mainardi, see \eg
\cite{GorMai INDIA03,GorMai CARRY04},
 have,
under proper  conditions on the tails
of the probability distributions for jumps
and waiting times,  investigated for Eq. (7.2)
the so-called  {\it well scaled transition to the diffusion limit}.
This transition is obtained by making smaller
all jumps by a positive factor $h$ and
 all waiting times by a positive factor $\tau $
related to $h$ by a {\it proper scaling relation},
and then letting  $h$ and $\tau$ tend to zero.
%%%%%%%%%%\vfill\eject
%%%%%%%%%%%%
An alternative interpretation is that we look at %% the one and
the same process
%% with a discrete number of jumps occurring after finite times,
from far away and after long time,
so that spatial distances and time intervals of normal size
appear very small but always related through
  {\it proper scaling}.
%%%%%%%%
In this limit the integral equation (7.2) has been shown to reduce
to a partial differential equation with fractional
derivatives in space and/or in time,
referred to as a {\it space-time fractional diffusion equation}.
This is also the topic
 of the recent paper by Scalas et al. \cite{Scalas PRE04} %% PRE03}
and, in a  more general framework, of the paper
by Gorenflo and Abdel-Rehim \cite{GAR Vietnam03}.
For a variant of well-scaled transition (there not called so) see
\cite{Uchaikin-Saenko 03}.
%%%%%%%%%%
In this paper  we limit ourselves to recall the
well scaled transition from the CTRW integral equation (7.2)
to the
{\it time fractional diffusion equation} (TFDE),
$$  %%% {\d^{\beta} \over \dt^{\beta} }\
  \null_tD_*^\beta \, u (x,t) :=
  {\ds \rec{\Gamma(1-\beta )}} \,
  {\ds \int_0^t}   \l[{\ds{\d \over \d \tau}} \, u (x,\tau)\r]\,
  {\ds{d\tau \over (t-\tau)^{\beta}}} ={\d^2 \over \d x^2} \,u (x,t)
 \,,
   \eqno(8.1)$$
with  $0<\beta \le  1$, subjected to the initial condition
$u (0,t)= \delta (x)$.
%% where  $\,_tD_*^\beta $ denotes the Caputo fractional derivative
%% of order $\beta$,see \eg \cite{GorMai CISM97,Podlubny BOOK99}.
%%%%%%%%%%%%
For  $\beta = 1$
we recover the standard diffusion equation, for which
the fundamental solution
is the Gaussian probability density evolving
in time with variance $\sigma ^2 = 2t$,
$$  u (x,t)
 = {1\over 2\sqrt{\pi }}\,t^{-1/2}\, \e^{-\ds x^2/(4t)}\,.
                            \eqno(8.2)$$
This probability law  is known to  govern the classical Brownian model
for the phenomenon of normal diffusion.
%%%%%%%%%
Also for $0<\beta <1$
the  fundamental solution of t (8.1)
can be interpreted  as a  (spatially-symmetric)  density
evolving in time,  but
  exhibiting stretched exponential tails with a variance
 proportional to $ t^{\beta }\,, $  implying a
phenomenon of {\it slow anomalous diffusion}.
%%%%%%%%%%%
We have, see for details
\cite{GoLuMa 99,GoLuMa 00,Mainardi CHAOS96,%%
Mainardi CISM97,Mainardi LUMAPA01,Mainardi MELFI01},
$$u(x,t)
   =    {1\over 2 \, t^{\beta/2} }\, M_{\beta /2}(|x|/t^{\beta /2})
\,,\eqno(8.3)$$
 where
$$M_{\beta/2}(y):=
   {\ds \sum_{n=0}^{\infty}{(-y)^n\over n!\,\Gamma[-\beta
     n/2+(1-\beta/2 )]}}\,,
   \eqno(8.4)$$
$|x|/t^{\beta /2} $ acting as a {\it similarity variable}.
%%%%%%%%%%%%\vfill\eject
%%%%%%%%%%%%%\noindent
We stress that, according to footnote $^{(4)}$,
 $M_{\beta/2}(y)=  \Phi_{-\beta /2,1-\beta /2}(y)$,
hence the fundamental solution of the TFDE
  is a  Wright function\footnote{
%%%%%%%%%%%%%%%%%%%%%%%%%%%%%%%%%%%%%%%%%%%%%%%%%%%%%%%%%%%%%%%%%%%%%%%%%%%%%%
In his  first 1993 analysis of the time fractional diffusion equation,
one of the present authors (F.M.) \cite{Mainardi WASCOM93}
  introduced the two
{\it auxiliary functions} (of Wright type),
$ F_\nu (z) :=   \Phi _{-\nu , 0}(-z)$ and
$ M_\nu (z) :=  \Phi _{-\nu , 1-\nu }(-z)$ with $ 0<\nu <1\,,$
inter-related through $ F_\nu (z) = \nu  \, z \, M_\nu (z ) \,.$
Being in that time only aware of the Bateman project
where the parameter $\lambda  $ of the Wright function
$\Phi_{\lambda,\mu}(z)$
was erroneously restricted to  non-negative values,
%% (in order to obtain an entire function)
F.M. thought to have extended the original Wright  function.
It was just Professor Stankovi{\'c}
during the presentation of the paper
\cite{MainaTomi 95}
in the Conference {\it Transform Methods and Special Functions, Sofia 1994},
who informed F.M. that this extension for $-1<\mu < 0$
was already made just by Wright himself in 1940 (following
his previous papers in 1930's).
%% \\
On this special occasion F.M. wants
to renew  his personal gratitude to Professor
Stankovi{\'c}  for this earlier information, that has induced
 him to study the original papers by Wright and work  (also in collaboration)
on the functions of the Wright type for further applications.}.
%%%%%%%%
The TFDE (8.1) can be derived from the CTRW integral equation (7.2),
by %%  compressing and
properly rescaling the  waiting time and the jump widths
and passing to the diffusion limit.
%%%%%%
Making smaller all waiting times by a
positive factor $\tau \,, $   all jumps
by a positive factor $h\,,$
we get, for $n\in \NN$,  the jump instants
$ t_n(\tau ) = \tau T_1 + \tau T_2 + \dots + \tau T_n$
and the jump sums, $ x_0(h )=0$,
$ x_n(h) = h X_1 + h X_2 +  \dots + h X_n$.
The reduced waiting times $\tau T_n$ all have the $pdf$
$\phi _\tau (t) = \phi (t/\tau )/\tau \,, \, t>0\,, $
 analogously the reduced jumps $hX_n$  all have
the $pdf$ $w_h(x) = w(x/h)/h \,, \, x \in \RR\,. $
%% Under this scaling transformation
The  probability density $p(x,t)$ so obtained
is denoted by $p_{h,\tau} (x,t)$.
Then, the transition to the diffusion limit consists
in sending $\tau \to 0$ and $h\to 0$ under an appropriate
relation between $\tau $ and $h$ and deriving the evolution
equation satisfied by  the limiting $pdf$  $p_{0,0}(x,t)$.

Let us now  resume the approach by Gorenflo and Mainardi,
for giving  the conditions on  $w(x)$, $\phi (t)$ and
the {\it scaling relation} between $h$ and $\tau $,
which allow   $p_{0,0}(x,t)$ to be identified  with
the fundamental solution $u(x,t)$ of our TFDE (8.1).
%%%%%%%%%%%%%%%%
Assuming the jump $pdf$ $w(x)$ to be symmetric   with  finite variance
$  \sigma ^2$
%% := \int_{-\infty}^{+\infty}  x^2\, w(x)\, dx$
and  the waiting time $pdf$
 $\phi(t) $  either with finite mean
$ \rho$
%% := \int_0^\infty t\, \phi(t)\, dt$
(relevant in the case $\beta =1$), or,
with $c>0 $ and some $\beta \in (0,1)\,, $
$ \phi(t) \sim c\, t^{-(\beta +1)}$
for  $t \to \infty$,
and setting
 $  \mu  = \sigma ^2$,
 $\lambda  = \rho$ if $\beta  =1$,
 $\lambda=  c\, \Gamma(1-\beta) / \beta$
   if $ 0<\beta  <1$,
the  required {\it scaling relation}
  for the diffusion limit is
$   \lambda \,\tau^{\beta} = \mu \, h^2 $.

We  now will compare, at fixed times, the    spatial
probability distributions  %%%  evolving in time
provided by the  fundamental solutions
of the integral equation of the CTRW as given by the series (7.5)
with the fundamental solutions of the corresponding limiting
TFDE as   given by the function (8.3) of Wright type.
%%%%%%%%%%\vfill \eject
%%%%%%%%%%%%%
This will enable us, for the compound Poisson, Mittag-Leffler
and Wright processes,
to investigate from a numerical view  the increasing quality of
approximation to the diffusion limit with advancing time,
independently from the scaling relation.
%%%%%%%%%%
In this numerical comparison, in order to avoid the singular terms
originated by delta functions,
 we prefer to consider,
 instead of the spatial densities $p(x,t)$, $w(x)$
and $u(x,t)$,
the corresponding  %%{\it probability distribution functions},
 {\it probability cumulative  functions}   %% ($pcf$)
that we denote by  capital letters, namely
$P(x,t)$, $W(x)$ and $U(x,t)$.
%%%%%%%%%%%%%%
Thus,  by integrating in space %% the terms in
(7.5) from $-\infty$ to $x$,
and denoting by
$\Theta(x)$ the unit step Heaviside function,
we have
%% the $pcf$ $P(x,t)$ is evolving in time according to
$$ P(x,t)  =  \Psi(t)\, \Theta(x) + \sum_{k=1}^{\infty} v_k(t)\, W_k (x)\,.
\eqno(8.5)$$
In an analogue way, for  the TFDE  we have
$$ U(x,t) =
 {1\over 2}\,\l[ 1 +
\int_{0}^{x/t^{\beta /2}}\!\!\!   M_{\beta /2}(y) \, dy\r]\,,
\q 0<\beta \le 1\,,
\eqno(8.6)$$
that in the case $\beta =1$ reduces to
$$ U(x,t)  = {\rec 2} \l[ 1 + {\erf} \l( {x\over 2\,\sqrt{t}}\r)\r]\,.
                            \eqno(8.7)$$
	%%%%%%%%%% FIGURES 4 5 6 %%%
\begin{figure}%% [h!]
 \includegraphics[width=.52\textwidth]{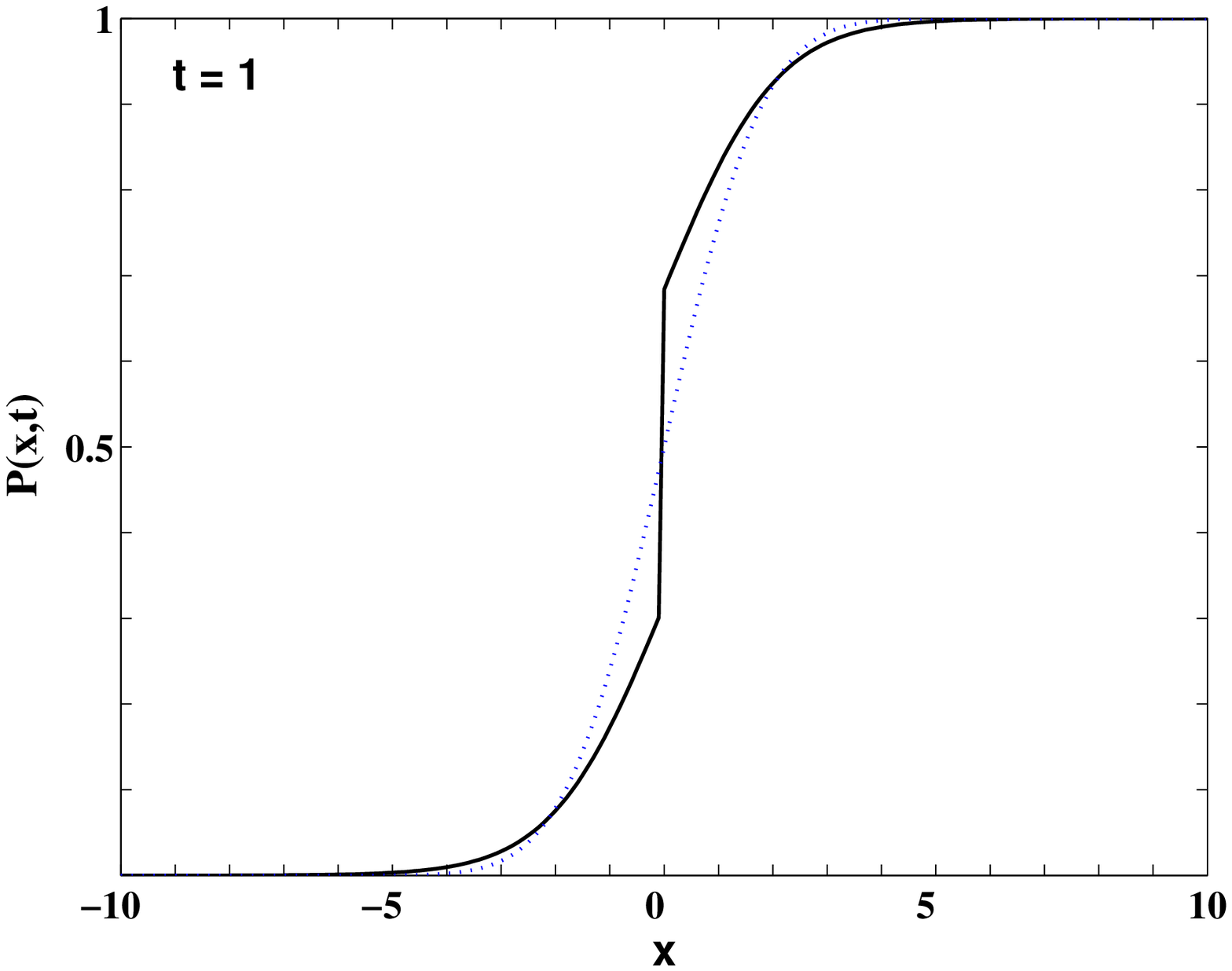}
\includegraphics[width=.52\textwidth]{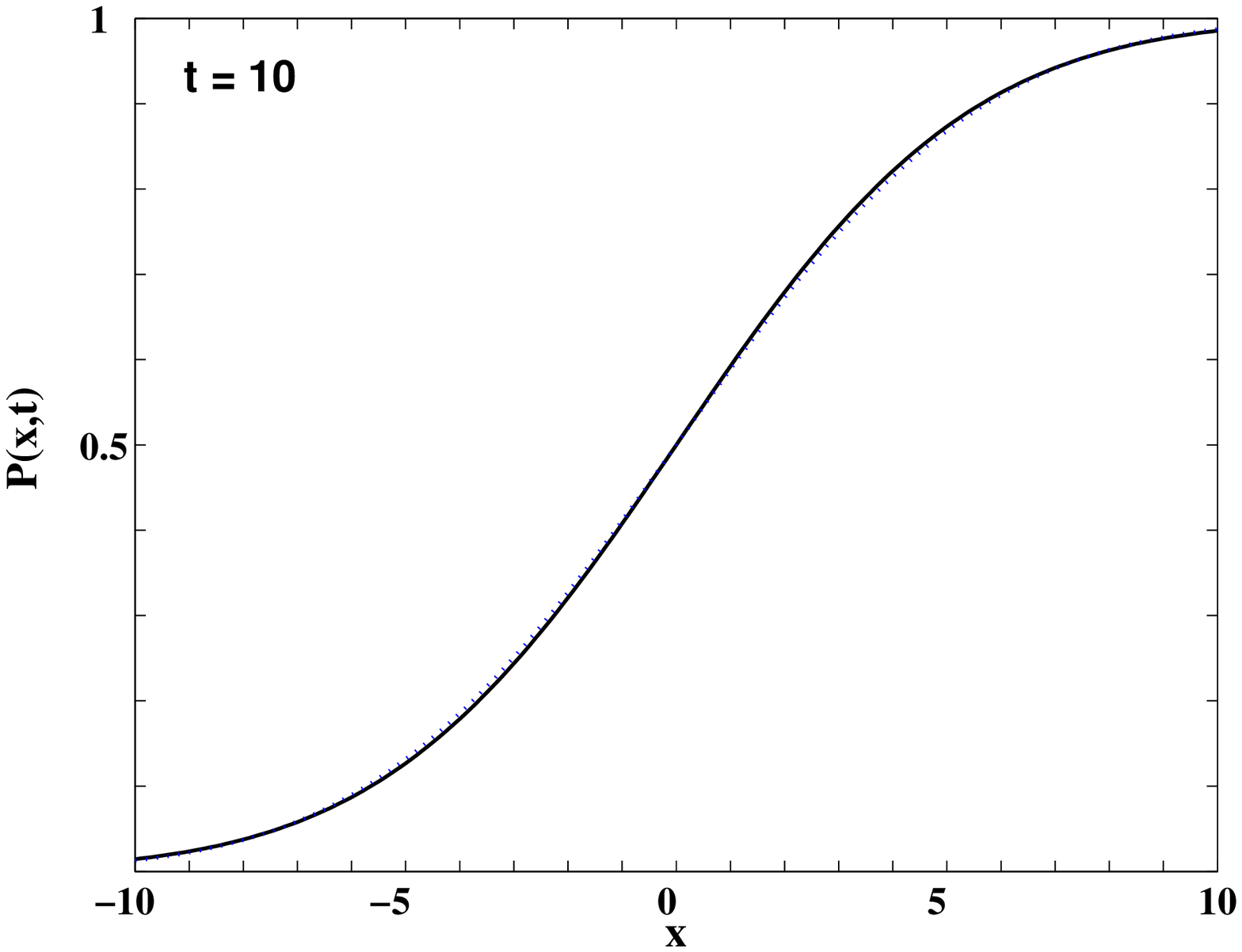}
\caption{Compound Poisson process.}
% \end{figure}
%\centerline{Left: $\Psi(t=1)=3.68\,10^{-1}$; Right: $\Psi(t=10)=4.54\,10^{-5}+
%%%%%%%%%%%%%%%%%%%%%%%%%%%%%%%%%%%%%%%%%%%%%%%%%%%%%%%%%%%%%%%%%%%%%%%
% \begin{figure}%% [htbp]
 \includegraphics[width=.52\textwidth]{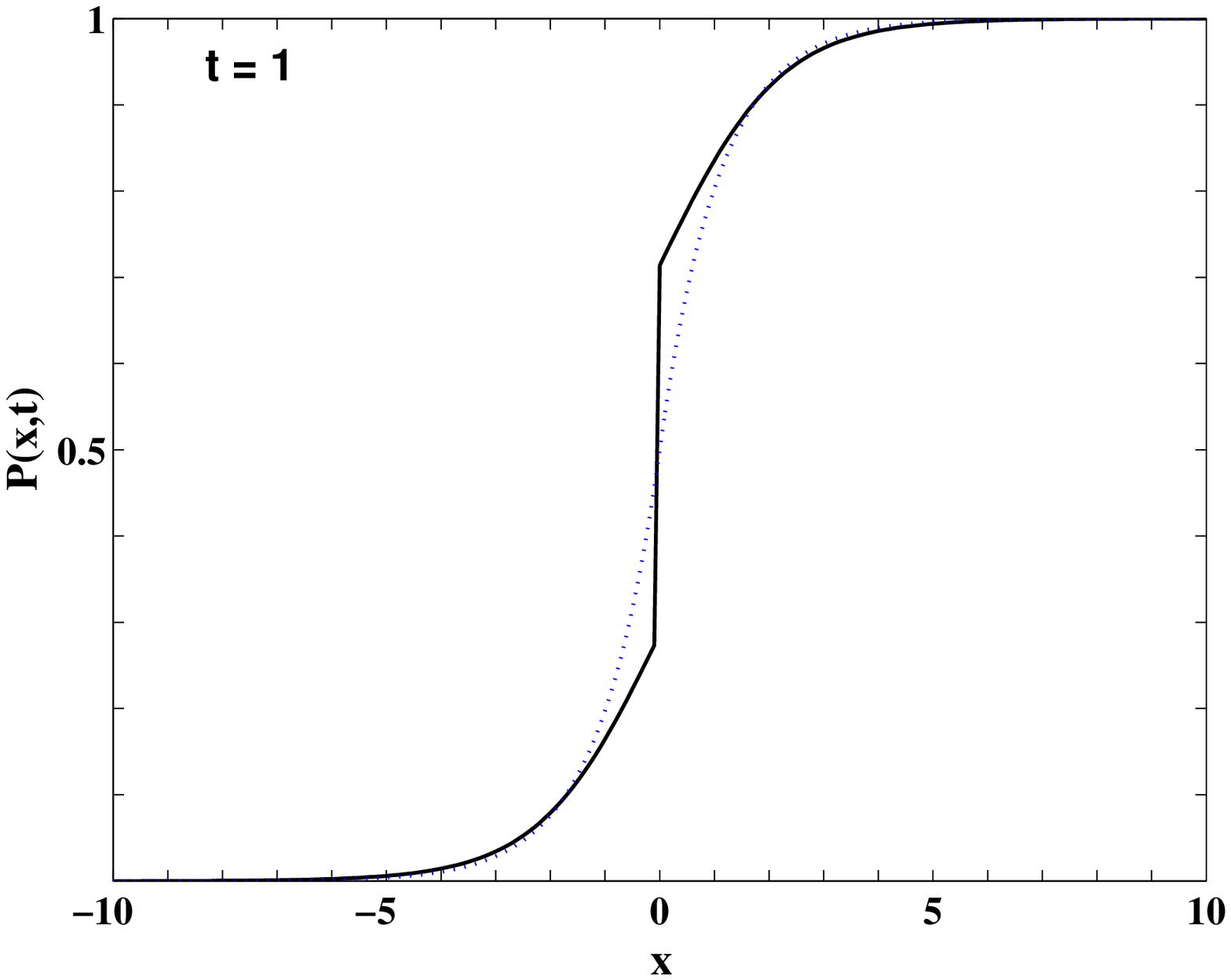}
\includegraphics[width=.52\textwidth]{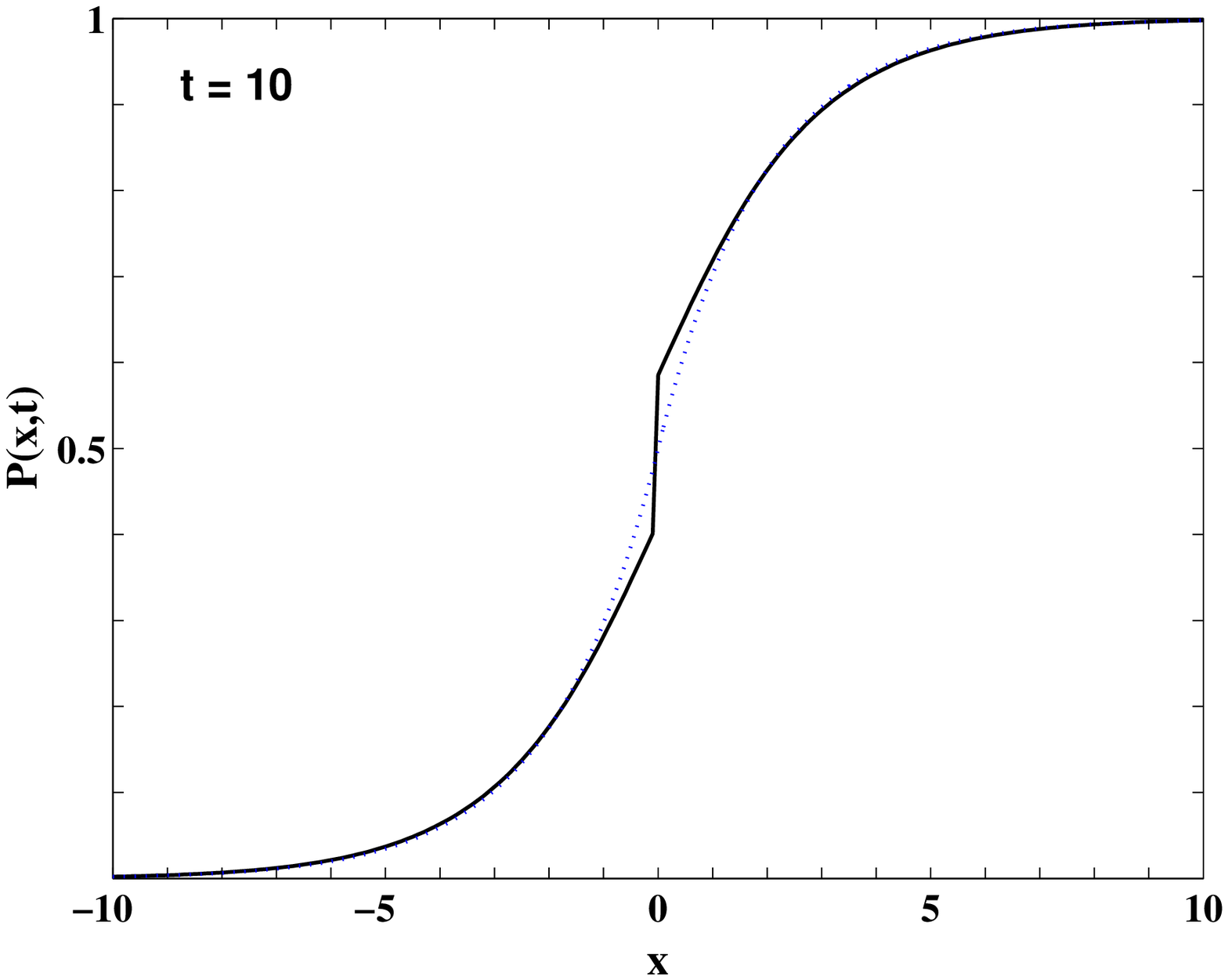}
\caption{Compound Mittag-Leffler  process.}
% \end{figure}
%\centerline{Left: $\Psi(t=1)=4.28\,10^{-1}$; Right: $\Psi(t=10)=1.71\,10^{-1}$
%%%%%%%%%%%%%%%%%%%%%%%%%%%%%%%%%%%%%%%%%%%%%%%%%%%%%%%%%%%%%%%%%%%
%\begin{figure}%% [htbp]
 \includegraphics[width=.52\textwidth]{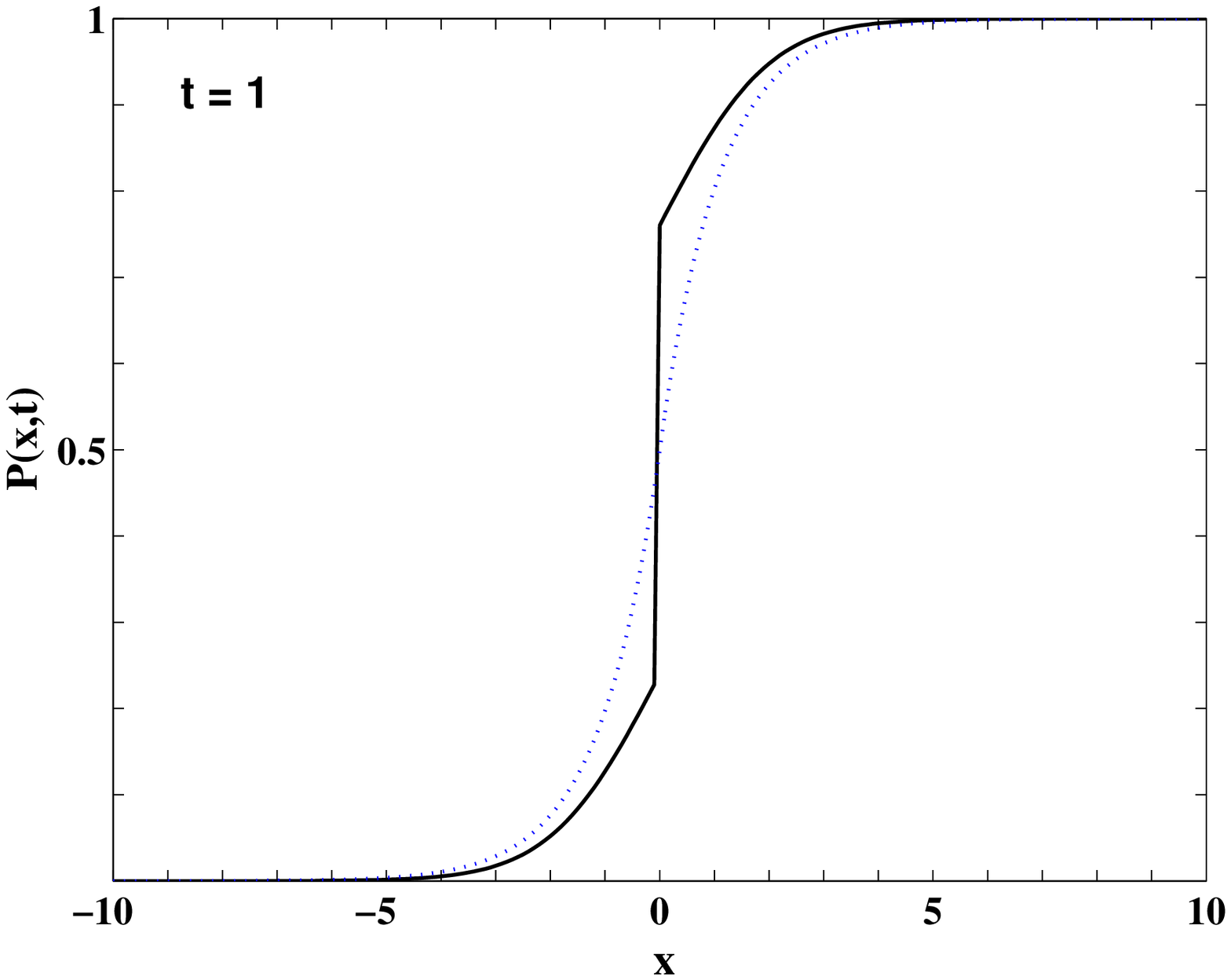}
\includegraphics[width=.52\textwidth]{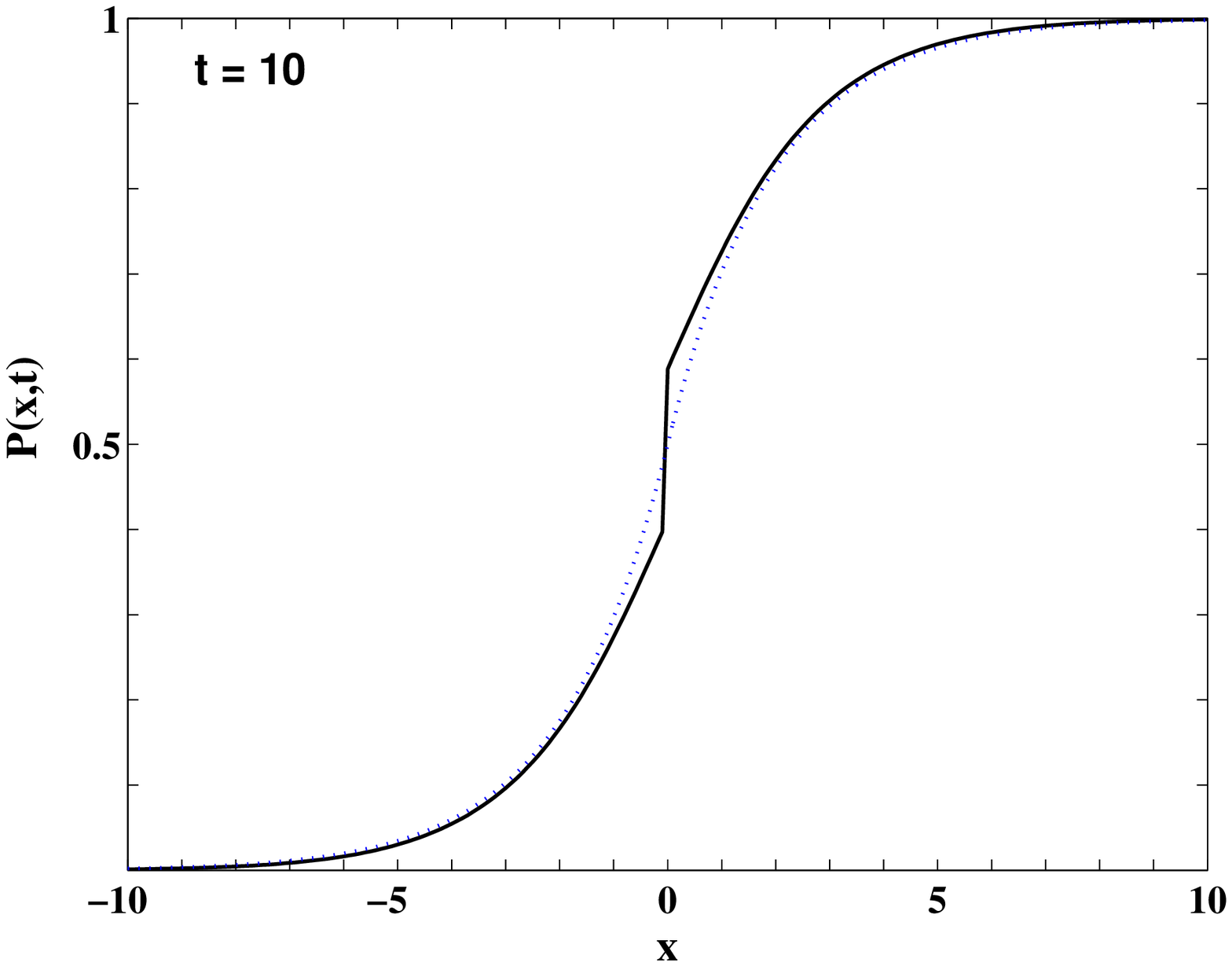}
\caption{Compound Wright  process.}
%\centerline{Left: $\Psi(t=1)=5.21\,10^{-1}$; Right: $\Psi(t=10)=1.77\,10^{-1}$
\end{figure}
%%%%%%%%%%%%% THE END OF FIGURES 4 5 6						
							
In our case studies we have
chosen $t=1$ and $t=10$,
 as typical short and long times.
As waiting time distributions
we have used the
exponential  distribution  for the compound Poisson process
so that the $v_k(t)$ are given by (3.5),
and the Mittag-Leffler  and   Wright distributions
for $\beta =1/2$, for the corresponding compound processes
so that $v_k(t)$ are given by (4.10) with (6.7), and (6.8), respectively.
For all the processes
 we have used  a {\it symmetric Gaussian}
for the common jump probability distribution,
for which we have
$$ w_k (x) =\frac{1}{2\sqrt{\pi}}\,
\frac{\e^{\, \ds -x^2/(4k)}}{\sqrt{k}} \,, \q
W_k (x) = \frac{1}{2}\,
 \l[ 1 + {\erf} \l( {x\over 2\,\sqrt{k}}\r)\r]\,.
                            \eqno(8.8)$$
%%%%%%%%%%%%%%%%

%%%%%%%%%%

%%%%%%%%%%%
In figures 4,5,6, regarding respectively,
the compound Poisson, Mittag-Leffler
and Wright processes,  we  report   in continuous line
 the $pcf$
for the CTRW whereas in dashed line
the corresponding cumulative functions  of the limiting TFDE.
%%%%%%%
The  values of the survival function $\Psi(t)$
for $t=1$ and $t=10$
are given
in Table I; however, they are
clearly visible as the height of the vertical part (at $x=0$)
of the continuous line.
%%%%%%%%%%%

We note that for $\beta =1/2$ there is no substantial difference
 between the cases of Mittag-Leffler and Wright compound processes:
 at these times the dominant effect
is due to the common  asymptotic behaviour of the corresponding distributions
resulting from the diffusion limit.
In other words:  %%%  we note that the
the  results are consistent with the asymptotic
properties of the diffusion limit.
We also note that in this comparison we cannot
take  direct profit of the scaling relation of the diffusion limit,
because the fundamental solution of the CTRW integral equation
has no self-similarity property at variance
to the fundamental solution of the corresponding   TFDE.
%% \vfill\eject

We  point out that also Barkai \cite{Barkai ChemPhys02}
has investigated the diffusion limit for the CTRW process of Wright type
but considering different case studies
and using different notation and terminology.

%%%%%%%%%%%%%%
\section*{9. Conclusions}

After sketching the basic principles of renewal theory we
have discussed the classical Poisson process.
It is well known that in this process the sum of $k$
%% ,$ k\in \NN$,
waiting times obeys the so-called {\it Erlang distribution},
named after the Danish telecommunication engineer Agner Krarup
{\it Erlang} (1878-1929)\footnote{
%%%%%%%%%%%%%%%
See
{\tt http://plus.maths.org/issue2/erlang/}
for information on {Erlang}.}.
%%%%%%%%%%%%%%%%%%%%%%%%%%%%%%%%

Because the exponential function is just a
special case of the
{\it Mittag-Leffler function}, one may expect that an analogous explicit
representation exists for the situation of waiting times distributed
according to the Mittag-Leffler law  $\Psi(t)=E_{\beta}(-t^\beta)$,
$0<\beta \le 1$, with   $\Psi(t)$ denoting the survival function (the
probability of no event occurring in an interval of length $t$ just
after an event). We have thus presented such explicit representations,
for the above renewal process that we have called
{\it renewal process of Mittag-Leffler type}.
%%%
We have also considered %%  offer
another process where
the survival function, being related to the unilateral
L\'evy stable law of index $\beta $,
is given as an
expression involving the {\it Wright function}.
This process has been called {\it renewal process of Wright type}.
%% namely the renewal process of Mittag-Leffler type and the renewal process
%% of Wright type, so named by us because special functions of
%% Mittag-Leffler and of Wright type appear in the definition of the relevant
%% waiting times.

In both of these renewal processes the waiting time probability distribution
is long-tailed, and there is a parameter $\beta$, $0 < \beta < 1$.
The first of these processes is a direct generalization of the Poisson
process as for the limiting  value $\beta = 1$ the Poisson process is
recovered. In distinct contrast, the second process changes its character
from stochastic to deterministic: for $\beta = 1$ at every instant $t=n$,
$n$ a natural number, an event happens (and never for other instants).
%%%%%%%%
Then,  we have reviewed the general theory of renewal processes with reward,
the so called {\it compound renewal processes}, known  in  physics
and chemistry literature  as {\it continuous time random walks}.

Using the fact that the sojourn probability density
%% (evolving in time)
of a compound renewal process
 can be written as an infinite series involving the
convolution powers of the waiting times and the rewards (the jumps),
we have simulated these processes (all three cases: first {\it Poisson},
then, with  $\beta=1/2$,
 {\it Mittag-Leffler} and  {\it Wright}).
%%%%%%%%
This comparison has allowed us to make visible how in the large time regime
the behaviour of the probability cumulative  functions approaches
that of the corresponding functions for
{\it time-fractional diffusion} processes, to which
our compound renewal processes are known to weakly converge
in the well scaled transition to the diffusion limit.
Nicely  one can see that the survival function decays
towards zero very fast (exponentially) in the Poisson process, but
slowly (like a power with negative exponent) in the other two
processes.
%%%%%%%%%%%
Incidentally  we have also stressed,
from a theoretical view-point, %% laid open
the remarkable contrast
in the behaviour of the non-Markovian compound
processes of Mittag-Leffler and of Wright
type in the limit  $\beta = 1$, yielding the classical compound
Poisson process
(which is continuous in time and Markovian)  for the first type,
but a process discrete in time
(Markovian on the set of the integer time-values,
but non Markovian on the set of the real time-values) for the second.
%%%%%%%%%%%%%
\vskip -0.55 truecm
%%%%%%%%%%%%
 \section*{Acknowledgements}   %% {\bf Acknowledgements}

This work has been  carried out in the framework of a joint research
project for {\it Fractional Calculus Modelling}
({\tt www.fracalmo.org}), with partial support
of the Italian COFIN-MIUR project
%% for{\it Chaos, Nonlinearity and Anomalous Transport}
{\it Order and Chaos in Nonlinear Extended Systems:
 Coherent Structures, Weak Stochasticity and Anomalous Transport}
({\tt www.physycom.unibo.it}).
R.G. and F.M.  appreciate the support
of the EU ERASMUS-SOCRATES program for visits to Bologna and Berlin.
%% that, besides for teaching, were also useful for this research activity.
%%% \end{acknowledgements}

%%%%%%%%%%%%%%%%%%%%%%%%
\vfill\eject

%%%%%%%%%% REFERENCES %%%%%%%%


\begin{thebibliography}{99}

%%%%%%%%%%%%%  $A$ %%%%%%%%%%%%%%

 \bibitem{AS 65}
 {M. Abramowitz and I.A. Stegun},
  {\it  Handbook of Mathematical Functions}.
 New York, Dover  (1965).


%%%%%%%%%%%%%%%  $B$ %%%%%%%%%%%%%%%


\bibitem{Baeumer-Meerschaert FCAA01}
 B. Baeumer and M.M. Meerschaert,
 Stochastic solutions for fractional Cauchy problems.
 {\it Fractional Calculus and Applied Analysis}
 {\bf 4} (2001),  481-500.



% \bibitem{Barkai PRE01}
% E. B a r k a i,
% {Fractional Fokker-Planck equation, solution and application}.
% {\it Phys. Rev. E} {\bf 63} (2001), 046118/1-17.




\bibitem{Barkai ChemPhys02}
E. Barkai,
CTRW pathways to the fractional diffusion equation.
{\it Chemical Physics} {\bf 284} (2002),  13--27.

\bibitem{Barkai PRE00}
E. Barkai,  R. Metzler and  J. Klafter,
From continuous time random walk to fractional
Fokker-Planck equation.
  {\it Physical Review E} {\bf 61} (2000), 132-138.

\bibitem{Bazhlekova FCAA00}
E.G. Bazhlekova,
Subordination principle  for fractional evolution equations.
{\it Fractional Calculus and Applied Analysis} {\bf 3} (2000), 213-230.


%%  \bibitem{Beichelt-Fatti STOCHPROC02}
%%  F.E. B e i c h e l t and L.P. F a t t i,
%%  {\it Stochastic Processes and Their Applications}.
%%  London and New York, Taylor \& Francis  (2002). %%% [X]
%%% Personal copy bought at Copenhagen May 2002 565 DKK
%%%%%%%%%%%%%%%%%%%   $C$ %%%%%%%%%%%%%%
\bibitem{Caputo FCAA01}
 M. Caputo,
 Distributed order differential equations modelling
  dielectric induction and diffusion.
{\it Fractional Calculus and Applied Analysis}
{\bf 4} (2001), 421-442.


% \bibitem{Caputo-Mainardi RNC71}
%{M. C a p u t o and  F. M a i n a r d i},
%  Linear models of dissipation in  anelastic solids.
%  {\it Riv. Nuovo Cimento} (Ser. II) {\bf 1} (1971), 161--198.

% \bibitem{ChechkinGorenfloSokolov PRE02}
% A.V. C h e c h k i n, R. G o r e n f l o and I.M. S o k o l o v,
% Retarding subdiffusion and accelerating superdiffusion
% governed by distributed-order fractional diffusion equations.
 {\it Phys. Rev. E} {\bf  66} (2002) 046129/1-6.

\bibitem{Chechkin FCAA03}
A.V. Chechkin, R. Gorenflo,  I.M. Sokolov and
V.Yu. Gonchar,
Distributed order time fractional diffusion equation.
{\it Fractional Calculus and Applied Analysis} {\bf 6} (2003), 259-279.


% \bibitem{ChechkinKlafterSokolov EUROPHYSICS03}   %%%
% A.V. C h e c h k i n, J. K l a f t e r and I.M. S o k o l o v,
% Fractional Fokker-Planck equation for ultraslow kinetics,
% {\it Europhysics Lett.} {\bf  63} (2003) 326-332.



\bibitem{Cox RENEWAL67}
{D.R. Cox},
{\it Renewal Theory}. 2-nd Edn London, Methuen (1967). %% pp. 142.
% [1-st Edn 1962]

% \bibitem{Cramer 46}
% H. Cram\'er,
% {\it Mathematical Methods of Statistics},
% Princeton University Press, Princeton, N.J. (1946).

%%%%%%%%%%%%%%%%    $D$ %%%%%%%%%%%%%%%


% \bibitem{Doetsch 74}
%  G. Doetsch,
% {\it Introduction to the Theory and Applications of the
% Laplace Transformation} (Springer Verlag, Berlin, 1974).

% \bibitem{Djrbashian 66}     %%% Dzherbashian,
% {M.M. D j r b a s h i a n},
% {\it Integral Transforms and Representations of
% Functions in the Complex Plane}.
% Moscow, Nauka (1966). In Russian.
% Note that there is also the transliteration as  Dzherbashian.

%%%%%%%%%%%%%%       $E$ %%%%%%%%%%%%

\bibitem{Embrechts 01}
{P. Embrechts, C. Kl\"uppelberg and T. Mikosch},
{\it Modelling Extreme Events for Insurance and Finance}.
Berlin, Springer Verlag  (2001). %% , pp. 647.
%% [Springer Series in Applied Mathematics, Stochastic Modelling and Applied
%% Probability, No. 33]
%%%%%%%%%%%%%%

\bibitem{Erdelyi HTF}
%%  A. Erd\'elyi (Editor)       Bateman Project,
A. Erd\'elyi, W. Magnus, F. Oberhettinger
 and F.G. Tricomi,
  {\it Higher Transcendental Functions},  Vol. 3.
New York, McGraw-Hill  (1953-1954).


%%%%%%%%%%        $F$ %%%%%%%%%%

% \bibitem{Feller RENEWAL41}
% W. Feller,
% On the integral equation of renewal theory,
% {\it Annals of Mathematical Statistics} {\bf 12} (1941) 243-267.

%%%%%%%%%%%
\bibitem{Feller 71}
{W. Feller},
{\it An Introduction to Probability Theory and its Applications}, Vol. 2,
 2-nd edn. New York, Wiley (1971). [1-st edn. 1966]



%%%%%%%%%%%%      $G$ %%%%%%%%


 \bibitem{Stankovic 76}
 Lj.  Gaji{\'c} and  B. Stankovi{\'c},
   Some properties of Wright's function.
   {\it Publ. de l'Institut Math\`ematique, Beograd, Nouvelle S\`er.}
  {\bf 20} (1976), 91-98.


%%  \bibitem{Gelfand-Shilov 64}
% {I.M. G e l' f a n d and G.E. S h i l o v},
% {\it Generalized Functions}, Vol. 1.
% New York, Academic Press  (1964).
%[Translated from the Russian, Moscow (1958)]

%% \bibitem{Gnedenko-Kolmogorov LIMITDISTR54}
%% B.V. G n e d e n k o and A.N. K o l m o g o r o v,
%% {\it Limit Distributions for Sums of Independently Random Variables}.
%% Cambridge (Mass), Addison-Wesley  (1954).
%% [English translation from the Russian edition,  G.I.T.T.L., Moscow,  1949.]



% \bibitem{Gnedenko-Kovalenko QUEUEING68}
% B.V. G n e d e n k o and I.N. K o v a l e n k o,
%  {\it Introduction to Queueing Theory}.
% Jerusalem,    Israel Program for Scientific Translations
%  (1968). %%%  pp. 281.
%  [Translated from the 1966 Russian edition]


% \bibitem{Gomes-Pestana 78}
% M.I. Gomes and D.D. Pestana,
% The use of fractional calculus in probability,
% {\it Portugaliae Mathematica} {\bf 37} No 3-4 (1978) 259--271. [X]


\bibitem{GAR Vietnam03}
{R. Gorenflo and E. Abdel-Rehim},
From power laws to fractional diffusion: the direct way.
{\it Vietnam Journal of Mathematics} {\bf 32} SI (2004), 65-75.

% \bibitem{GoIsLu 00}
% {R. G o r e n f l o, A. I s k e n d e r o v and   Yu. L u c h k o},
% Mapping between solutions of fractional diffusion-wave equations.
%{\it Fractional Calculus \& Applied Analysis\/}
%    {\bf  3}  (2000), 75--86.


\bibitem{GoLuMa 99}
{R. Gorenflo, Yu. Luchko and  F. Mainardi},
Analytical properties and applications of the Wright function.
{\it Fractional Calculus and  Applied Analysis} {\bf 2} (1999), 383-414.

\bibitem{GoLuMa 00}
{R. Gorenflo, Yu. Luchko and   F. Mainardi},
Wright functions as scale-invariant solutions of the diffusion-wave
 equation.
{\it J. Computational and Applied Mathematics}
     {\bf 118}  (2000), 175-191.


%\bibitem{GoLuRo 97}
% {R. G o r e n f l o, Yu. L u c h k o and  S. R o g o s i n},
% Mittag-Leffler type func\-tions: notes on growth properties
%and distribution of zeros,
%{\it   Pre-print  A-04/97,
% Fachbereich Mathematik und Informatik, Freie Universit\"at, Berlin} (1997).
% \\ {\tt http://www.math.fu-berlin.de/publ/index.html}


\bibitem{GorMai CISM97}
{R. Gorenflo and F. Mainardi},
Fractional calculus:
integral and differential equations of
fractional order. In
A. Carpinteri and  F. Mainardi (Editors),
{\it Fractals and Fractional
Calculus in Continuum Mechanics}.
Wien and New York, Springer Verlag (1997), pp. 223-276.
%% [CISM Lecture Notes Vol. 378]. %% \\
%% Reprinted in
[{\tt http://www.fracalmo.org}]

\bibitem{GorMai CHEMNITZ01}
 {R. Gorenflo  and  F. Mainardi},
 Random walk models approximating symmetric
 space fractional diffusion  processes.
  In  J. Elschner, I.  Gohberg  and  B. Silbermann (Editors),
 {\it Problems in Mathematical  Physics}.
 Basel,   Birkh\"auser Verlag (2001), pp. 120-145.
% [Operator Theory: Advances and Applications, No 121]
%%% Chemnitz.tex paper:  Siegfried Pr\"ossdorf Memorial Volume,
%% Proceedings of the 11-th TMP Conference


% \bibitem{GorMai AARHUS02}
%{R. G o r e n f l o and F. M a i n a r d i},
%     Non-Markovian random walk models, scaling and diffusion limits.
%In O. E. Barndorff-Nielsen (Editor),
%     Mini Proceedings of the 2-nd %% MaPhySto
%    Conference
%     on {\it L\'evy Processes: Theory and Applications}.
%     MaPhySto Centre, %%% (Mathematical Physics and Stochastics Centre),
%     %% Dept. Mathematics,
%    University of Aarhus (Denmark),  21-25 January 2002,
%     Miscellanea No. 22 (2002), pp. 120-128.
%    See {\tt http://www.maphysto.dk}


 \bibitem{GorMai INDIA03}
{R.  Gorenflo  and  F. Mainardi},
 Fractional diffusion processes: probability distributions and
    continuous time random walk. In
   G. Rangarajan and  M. Ding (Editors),
   {\it Processes with Long Range Correlations}.
 Berlin,  Springer Verlag  (2003), pp. 148-166.
 %% [Lecture Notes in Physics, No. 621]


%% \bibitem{GorMai FDA04}
%% {R. G o r e n f l o and F. M a i n a r d i},
%% Power laws, random walks, and fractional diffusion processes
%%   as well scaled  refinement limits.
%% In  A. Le M\'ehaut\'e, J.A. Tenreiro Machado, J.C. Trigeassou, J. Sabatier
%% (Editors),
%% Proceedings of the {\it First IFAC Workshop on
%% Fractional Differentiation and its
%% Applications (FDA'04)}. Bordeaux (France) 19-21 July 2004, pp. 36-47.
%% [Plenary lecture]

\bibitem{GorMai CARRY04}
{R. Gorenflo and  F.  Mainardi},
Simply and multiply scaled
diffusion limits for continuous time random walks.
{\it IOP (Institute of Physics) Journal of Physics: Conference Series},
   {\bf 7} (2005), 1-16
[S. Benkadda, X. Leoncini and  G. Zaslavsky (Editors),
{Proceedings of the  International Workshop on Chaotic Transport and
   Complexity  in Fluids and Plasmas},
   Carry Le Rouet (France), 20-25 June 2004]
   


\bibitem{Gorenflo 01}
R. Gorenflo, F.  Mainardi,
E. Scalas and M. Raberto,
{Fractional calculus and continuous-time finance III:
the diffusion  limit}. In  M. Kohlmann and S. Tang (Editors),
{\it Trends in Mathematics} - {\it Mathematical Finance}.
Basel,   Birkh\"auser  (2001), pp. 171-180.

%%%%%%%%%%%%%%%%%     $H$ %%%%%%%%%%

\bibitem{Hilfer PhysA03}    %%
 R. Hilfer,
 On fractional diffusion and continuous time  random walks.
 {\it Physica A}  {\bf 329} (2003) 35-39.

\bibitem{HilferAnton 95}
{R. Hilfer and L. Anton},
Fractional master equations and fractal time random walks.
{\it Phys. Rev. E} {\bf 51} (1995), R848--R851.


% \bibitem{Huillet 02a}
% T. Huillet,
% {Renewal processes and the Hurst effect},
% {\it J. Phys. A} {\bf  35} (2002), 4395-4413.

% \bibitem{Huillet 02b}
% T. Huillet,
% {On the waiting time paradox and related topics},
% {\it Fractals} {\bf 10} (2002), 1-20.

%%%%%%%%%%%%%%%%%%%%%    $I$  %%%%%%%%%%


%%%%%%%%%%%%%%%%%%%%%%%%%   $J$  %%%%%%%%%%%%%%


%%%%%%%%%%%%%%%%%  $K$ %%%%%%%%%%%%%%%%%%
% \bibitem{Khintchine QUEUEING60}
%  A.Ya. K h i n t c h i n e,
%  {\it Mathematical Methods in the Theory of Queueing}.
%  London,  Charles Griffin  (1960). %% pp. 120.
% [Translated from the 1955 Russian edition]


% \bibitem{Khintchine LIMITDISTR38}
% A.Ya. Khintchine, %% (1938)
% {\it Limit Distributions for the Sum of Independent  Random Variables}.
% O.N.T.I., Moscow, 1938, pp. 115.
% [in Russian: the English translation from the Russian
% by S. Rogosin is now available]  %%  [X]

% \bibitem{Kilbas-Saigo 96}
% {A.A. K i l b a s and M. S a i g o},
% On Mittag-Leffler type functions,
%     fractional calculus operators and solution of integral equations.
% {\it Integral Transforms and Special Functions}
% {\bf 4} (1996), 355-370.

\bibitem{KilbasSaigoTrujillo 02}
{A.A. Kilbas, M. Saigo and  J.J. Trujillo},
On the generalized Wright function.
{\it Fractional Calculus and  Applied Analysis} {\bf 5} %% No. 4
(2002), 437-460.



\bibitem{Kiryakova 94}
{V. Kiryakova},
   {\it Generalized Fractional Calculus and Applications}.
      Harlow, Longman (1994).
[Pitman Research Notes in Mathematics, Vol. 301]



 %%%%%%%%%%%  $L$ %%%%%%%%



%%%%%%%%%%%%%%%   $M$ %%%%%%%%%%%%%%%%
\bibitem{Mainardi WASCOM93}
{F. Mainardi},
 On the initial value problem for the fractional diffusion-wave equation.
 In S. Rionero and T. Ruggeri (Editors),
{\it Waves and Stability in Continuous Media}.
 Singapore,  World Scientific  (1994), pp. 246-251.
%%  [Proc. VII-th  WASCOM, Int. Conf. "Waves and Stability in Continuous
%%   Media",  Bologna, Italy, 4-7 October 1993]


% \bibitem{Mainardi AML96}
% {F. M a i n a r d i},
% The fundamental solutions for the fractional diffusion-wave equation.
% {\em Applied Mathematics Letters\/} {\bf 9} No 6 (1996), 23--28.

\bibitem{Mainardi CHAOS96}
{F. Mainardi},
  Fractional relaxation-oscillation and fractional
  diffusion-wave phenomena.
  {\it Chaos, Solitons and Fractals} {\bf 7} (1996), 1461--1477.

\bibitem{Mainardi CISM97}
{F. Mainardi},
Fractional calculus: some basic problems in
continuum and statistical mechanics.
In A. Carpinteri and  F. Mainardi (Editors),
{\it Fractals and Fractional
Calculus in Continuum Mechanics}.
Wien and New York, Springer Verlag (1997), pp. 291-348.
%% [CISM Lecture Notes Vol. 378]
%%  Reprinted in
[{\tt http://www.fracalmo.org}]

\bibitem{MaiGor JCAM00}
{F. Mainardi and R. Gorenflo},
On Mittag-Leffler-type functions in fractional evolution processes.
{\it  J. Computational and Applied Mathematics}
{\bf 118}  (2000), 283-299.

%%%%
%\bibitem{Mainardi FRACTAL04}
% {F. M a i n a r d i, R. G o r e n f l o and E. S c a l a s},
%    A renewal process of Mittag-Leffler type.
% In M. Novak (Editor),
%    {\it Thinking in Patterns: Fractals and Related Phenomena in Nature}.
% Singapore,    World Scientific  (2004), pp. 35-46
   %%% (ISBN 981-238-822-2).

\bibitem{Mainardi VIETNAM04}
F. Mainardi, R. Gorenflo and  E. Scalas,
    A fractional generalization of the Poisson processes.
    {\it Vietnam Journal of Mathematics} {\bf 32} SI (2004), 53-64.

\bibitem{Mainardi LUMAPA01}
{F. Mainardi, Yu. Luchko and G. Pagnini},
 The fundamental solution of the space-time fractional diffusion equation.
   {\it Fractional Calculus and  Applied Analysis}
 {\bf 4}  (2001), 153-192.
%% [Reprinted in NEWS 010401, see {\tt http://www.fracalmo.org}]
%% [Paper dedicated to Professor Rudolf Gorenflo for his 70-th birthday]

\bibitem{Mainardi MELFI01}
{F. Mainardi and G. Pagnini},
The Wright functions as solutions of the time-fractional diffusion
      equations.
{\it Applied Mathematics and Computation} {\bf 141}  (2003), 51-62.
%%     Guest Editor: H.M.R Srivastava, G. Dattoli and P.E. Ricci,
%%     Chairmen of the International Conference on "Advanced Special Functions
%%     and Applications" Melfi, Italy, June 24-229, 2001.

\bibitem{Mainardi FCAA03}
{F. Mainardi,  G. Pagnini and R. Gorenflo},
Mellin transform and subordination laws  in fractional diffusion processes.
{\it Fractional Calculus and Applied Analysis} {\bf 6} (2003), 441-459.


% \bibitem{Mainardi OPSFA03}
%{F. M a i n a r d i, G. P a g n i n i and R.K. S a x e n a},
%    Fox $H$ functions in fractional diffusion,
%    {\it J. Computational \& Appl. Mathematics} {\bf 178} (2005) 321-333.

%%%%%%%%%%%%%%%%
 \bibitem{Mainardi PhysA00}
F. Mainardi, M. Raberto, R. Gorenflo and
E. Scalas,
 Fractional calculus and continuous-time finance II: the waiting-time
     distribution.    {\it Physica A} {\bf 287}  (2000), 468--481.
%%%%%


\bibitem{MainaTomi 95}
 {F. Mainardi and M. Tomirotti},
  On a special function
  arising  in  the time fractional diffusion-wave equation.
  In P. Rusev, I. Dimovski and V. Kiryakova (Editors),
  {\it Transform Methods and Special Functions, Sofia 1994}.
  Singapore, Science Culture Technology  (1995), pp. 171--183.


%\bibitem{Mainardi FNL05}
%F. M a i n a r d i, A. V i v o l i and R. G o r e n f l o,
%  Continuous time random walk  and time fractional diffusion:
%     a numerical comparison between the fundamental solutions.
%{\it Fluctuation and Noise Letters} (2005), to appear.


%%  \bibitem{Mantegna 00}   %% [26] %%
%%     R.N. M a n t e g n a and H.E. S t a n l e y,
%%  {\it An Introduction to Econophysics: Correlations and Complexity in
%% Finance}.
%%  Cambridge,  Cambridge University Press (2000).

% \bibitem{Marichev 83}
% O.I. M a r i c h e v,
% {\it Handbook of Integral Transforms of Higher Transcendental Functions,
% Theory and Algorithmic Tables}.
% Chichester, Ellis Horwood (1983).



%\bibitem{Masoliver 03a}
% J. Masoliver, M. Montero and G.H. Weiss,
% {Continuous-time random-walk model for financial distributions},
% {\it Phys. Rev. E} {\bf 67} (2003) 021112/1--9.

\bibitem{M3 PRE02sub}
M.M. Meerschaert, D.A. Benson, H.P  Scheffler and
B. Baeumer,
Stochastic solutions of space-fractional diffusion equation.
{\it Phys. Rev. E} {\bf  65} (2002), 041103/1-4.
%%A PRESENTATION, GOOD AND SHORT,FOR SUBORDINATION


% \bibitem{M3 PRE02sol}
% M.M. M e e r s c h a e r t, D.A. B e n s o n, H.P  S c h e f f l e r and
% P. B e c k e r - K e r n,
% Governing equations and solutions of anomalous random walk limits.
% {\it Phys. Rev. E} {\bf  66} (2002) 060102-1/4.




\bibitem{MetzlerKlafter PhysRep00}
 R. Metzler and J. Klafter,
The random walk's guide to anomalous diffusion: a fractional dynamics
 approach.
{\it Physics Reports}  {\bf 339} (2000), 1-77.


 \bibitem{MetzlerKlafter JPhysA04}
  R. Metzler and J. Klafter,
  The restaurant at the end of the random walk: Recent developments
  in the description of anomalous transport by fractional dynamics.
 {\it J. Phys. A. Math. Gen.}  {\bf 37} (2004), R161-R208.
%%%%%%%%%%%%%%%%%%%%%%%%%%%%%%%%%%%%%%



% \bibitem{Metzler PRE98}
%  R. M e t z l e r, J. K l a f t e r, I.M. S o k o l o v,
% Anomalous transport in external fields:
% Continuous time random walks and fractional
% diffusion equations extended.
% {\em Physical Review E} {\bf 58} (1998), 1621-1633.




%% \bibitem{Miller-Samko CM97}
%%  K.S. Miller and S.G.  Samko,
%% A note on the Complete monotonicity of
%% the generalized Mittag-Leffler function,
%%  {\it Real. Anal. Exchange} {\bf 23}  (1997) 753-755.

%\bibitem{Miller-Samko CM01}
% {K.S. M i l l e r and S.G.  S a m k o},
% Completely monotonic functions,.
% {\it Integr. Transf. and Spec. Funct.} {\bf 12} (2001), 389-402.


 \bibitem{MontrollScher 73}
  E.W. Montroll and H. Scher,
 Random walks on lattices, IV:
 Continuous-time walks and influence of absorbing boundaries,
  {\it J. Stat. Phys.} {\bf 9} (1973), 101-135.




%% \bibitem{MontrollShles 84}
%%  E.W. Montroll and M.F. Shlesinger,
%%  On the wonderful world of random walks, in
%%  J. Leibowitz and E.W. Montroll (Editors),
%%  {\it Nonequilibrium Phenomena II: from Stochastics to Hydrodynamics},
%%  North-Holland, Amsterdam, 1984, pp. 1-121.
%%  [{Studies in Statistical Mechanics}, vol. XI]

 \bibitem{MontrollWeiss 65}
 {E.W. Montroll and G.H. Weiss},
  Random walks on lattices, II.
 {\it J. Math. Phys.}  {\bf 6} (1965), 167--181.

%% \bibitem{MontrollWest 79}
%%  E.W. Montroll and  B.J. West,
%%  On an enriched collection of stochastic processes, in
%%  E.W. Montroll and J. Leibowitz (Editors),
%%  {\it Fluctuation Phenomena},
%% North-Holland, Amsterdam, 1979,  pp. 61-175.
%%  [{Studies in Statistical Mechanics}, Vol. VII]

%%%%%%%%      $N$ %%%%%%%%%%

% \bibitem{Naber 04}
% M. N a b e r,
%  {Distributed order fractional subdiffusion}.
%  {\it Fractals} {\bf 12} (2004) 23-32.



%%%%%%%%%%%%% $O$ %%%%%%%%%

%\bibitem{Oberhettinger  TABLES73}
%  F. Oberhettinger,
%  {\it Fourier Transforms and Their Inverse},
%  Academic Press, New York and London  (1973).


%%%%%%%%%%%%%% $P$ %%%%%%%%%%%%%%%%%

% \bibitem{Pillai 90}
% R.N. Pillai,
% On Mittag-Leffler functions and related distributions,
% {\it Ann. Inst. Statist. Math.} {\bf 42} No. 1 (1990)  157--161. [X]




\bibitem{Podlubny BOOK99}
 {I. Podlubny},
  {\it Fractional Differential Equations}.
San Diego,  Academic Press  (1999).


% \bibitem{Prudnikov IS3}
% {A.P. P r u d n i k o v, Y.A. B r y c h k o v and O.I. M a r i c h e v},
% {\it Integrals and Series}, Vol. 3.
% New York, Gordon and Breach Science Publ. (1990).


%%%%%%%%%%%%%   $Q$ %%%%%%%%%%%%%%




%%%%%%%%%%%%%%%%   $R$ %%%%%%%%%%%%%%%%%%%%%

% \bibitem{Raberto 02}
% M. Raberto, E. Scalas, and F. Mainardi,
%  Waiting-times and returns in high-frequency financial data:
% an empirical study,
% {\it Physica A} {\bf 314} (2002) 749--755.


% \bibitem{Ross STOCHPROC96}
%  S.M Ross,
%   {\it Stochastic Processes}, 2-nd Edn,
%   Wiley, New York (1996).


\bibitem{Ross PROBMOD97}
 {S.M. Ross},
{\it Introduction to Probability Models},
6-th Edn.
New York, Academic Press (1997).

%%%%%%%%%%%%%%%%%    $S$ %%%%%%%%%%%%%%%%%%%%

\bibitem{Saichev 97}
 A.I.  Saichev  and  G.M. Zaslavsky,
 Fractional kinetic equations: solutions and applications.
 {\it Chaos} {\bf 7} (1997), 753--764.



\bibitem{SKM 93}
{S.G.  Samko, A.A. Kilbas and O.I. Marichev},
{\it Fractional Integrals and Derivatives: Theory  and  Applications}.
New York and London,  Gordon and Breach Science Publishers (1993).
%[Translation from the Russian edition,
%  Minsk, Nauka i Tekhnika (1987)]



%%%%%%%%%%%%%%%%%%

\bibitem{Scalas PhysA00}
{E. Scalas, R. Gorenflo and F. Mainardi},
{Fractional calculus and continuous-time finance}.
{\it Physica A} {\bf 284} (2000), 376--384.

%% \bibitem{Scalas FRACTALS03}
%% E. Scalas, R. Gorenflo, F. Mainardi and M. Raberto,
%% Revisiting the derivation of the fractional diffusion equation,
%% {\it Fractals} {\bf 11} Suppl. S (2003) 281-289.

\bibitem{Scalas PRE04}
{E. Scalas, R. Gorenflo and F. Mainardi},
 Uncoupled continuous-time random walks: solution
 and limiting behaviour of the master equation.
{\it Physical Review E} {\bf 69} (2004), 011107/1-8.

%\bibitem{Schneider 96}
% W.R. Schneider,
% Completely monotone generalized Mittag-Leffler functions,
% {\it Expositiones Mathematicae} {\bf 14} (1996) 3-16.

\bibitem{SchneiderWyss 89}
{W.R. Schneider and W. Wyss},
Fractional diffusion and wave equations.
 {\it J. Math. Phys.} {\bf 30} (1989), 134-144.




\bibitem{Sokolov PhysicaA01}
 I.M. Sokolov, A. Blumen and J. Klafter,
 Linear response in complex systems: CTRW
 and the fractional Fokker-Planck equations.
 {\it Physica A} {\bf 302} (2001), 268-278.


\bibitem{Sokolov 04}  %%%{SChK}
 I.M. Sokolov, A.V. Chechkin and J. Klafter,
Distributed-order fractional kinetics.
{\it Acta Physica Polonica} {\bf 35} (2004), 1323-1341.
% [{\tt http://arXiv.org/abs/cond-mat/0401146}]
% [E-print arXiv: {\tt cond-mat/0401146}] %% (v1, 9 Jan 2004, 18 pages)
%% Presented at the 16th Marian Smoluchowski Symposium
% on Statistical Physics: Fundamentals and Applications, September 6-11, 2003]


%\bibitem{Srivastava H}
% {H.M. S r i v a s t a v a, K.C. G u p t a and S.P. G o y a l}, %% \\
%{\it The H-Functions of One and Two Variables with Applications}.
%New Delhi and Madras, South Asian Publishers  (1982).

\bibitem{Stankovic 70}
{B. Stankovi\`c},
{On the function of E.M. Wright}.
{\it Publ. de l'Institut
Math\`ematique, Beograd, Nouvelle S\`er.} {\bf 10} (1970), 113-124.



%%%%%%%%%%%%%%%%%%%   $T$ %%%%%%%%%%%%%

% \bibitem{Thu 84}
% N. Thu, Fractional calculus in probability,
% {\it Probability and Math. Statistics} {\bf 3} No 2 (1984) 173-189. [X]

%%%%%%%%%%%%%%%%%%%%  $U$ %%%%%%%%%%%%%%%%%%%


\bibitem{Uchaikin-Saenko 03}
  V.V. Uchaikin and V.V. Saenko,
Stochastic solution of partial differential equations
of fractional order.
{\it Siberian  Journal of Numerical Mathematics}  {\bf 6} (2003), 197-203.

% \bibitem{Uchaikin-Zolotarev 99}
% V.V. Uchaikin, V.M. Zolotarev,
% {\it Chance and Stability: stable laws
% and their application},
%  VSP [Modern Probability and Statistics], Utrecht, 1999.

%%%%%%%%%%%%%%%%%%%        $W$ %%%%%%%%%%%%%%%%%%

 \bibitem{Weiss BOOK94}
 {G.H. Weiss},
  {\it Aspects and Applications of Random Walks}.
  Amsterdam,  North-Holland (1994).
%%%%%%

\bibitem{Wyss-Wyss FCAA01}
M.M. Wyss and W. Wyss,
Evolution, its fractional extension and generalization.
{\it Fractional Calculus and Applied Analysis} {\bf 4} (2001), 273-284.
\end{thebibliography}
\end{document}